%% file: __Abelian_Theory.tex
\documentclass[11pt,oneside]{book}
\usepackage{geometry}                		
\usepackage{graphicx}									
\usepackage{amssymb}\vspace{0.25cm}
\usepackage{mathtools}
\usepackage{setspace}
\usepackage{tikz-cd}
\usepackage[mathscr]{euscript}
\newcommand{\abs}[1]{\left\vert#1\right\vert}
\usepackage[normalem]{ulem}
\usepackage{manfnt}
\usepackage{pgfplots}
\usepackage{pifont}
\usepackage[safe]{tipa}
\usepackage{marvosym}
\usepackage{scalerel,stackengine}
\usepackage{titlesec}
\usepackage{makeidx}
\usepackage{centernot}

\usepackage{xspace}
\usepackage[noauto]{chappg}
\usepackage[OT2,T1]{fontenc}
\usepackage{hyperref}
\usepackage{amsmath}
\usepackage[amsmath,thmmarks]{ntheorem}

\newcommand\mLT{%
\raisebox{0.1cm}{$L$}$T$\xspace
}

\newcommand\LY{%
\raisebox{0.05cm}{\text{\textlhtlongy}}\xspace
}

\newcommand\mA{%
$A$\xspace
}
\newcommand\mB{%
$B$\xspace
}

\newcommand\mF{%
$F$\xspace
}
\newcommand\mG{%
$G$\xspace
}
\newcommand\mH{%
$H$\xspace
}

\newcommand\mM{%
$M$\xspace
}

\newcommand\mR{%
$R$\xspace
}
\newcommand\mS{%
$S$\xspace
}
\newcommand\mT{%
$T$\xspace
}
\newcommand\mU{%
$U$\xspace
}

\newcommand\mX{%
$X$\xspace
}
\newcommand\mY{%
$Y$\xspace
}

\newcommand\A{%
\mathbb{A}
}
\newcommand\C{%
\mathbb{C}
}

\newcommand\I{%
\mathbb{I}
}
\newcommand\K{%
\mathbb{K}
}
\newcommand\bk{%
\textbf{k}\xspace
}
\newcommand\LL{%
\mathbb{L}
}
\newcommand\N{%
\mathbb{N}
}

\newcommand\Q{%
\mathbb{Q}
}
\newcommand\R{%
\mathbb{R}
}

\newcommand\Z{%
\mathbb{Z}
}
\newcommand\pe{%
\mathfrak{p}\xspace
}
\newcommand\mfm{%
\mathfrak{m}\xspace
}
\newcommand\bAB{%
\textbf{AB}\xspace
}
\newcommand\bAFF{%
\textbf{AFF}\xspace
}
\newcommand\bAFFGRP{%
\textbf{AFF-GRP}\xspace
}
\newcommand\bALG{%
\textbf{ALG}\xspace
}
\newcommand\bET{%
\textbf{\'ET}\xspace
}
\newcommand\bGRP{%
\textbf{GRP}\xspace
}
\newcommand\bMOD{%
\textbf{MOD}\xspace
}
\newcommand\bRNG{%
\textbf{RNG}\xspace
}
\newcommand\bSCH{%
\textbf{SCH}\xspace
}

\newcommand\bSET{%
\textbf{SET}\xspace
}

\newcommand\bTOR{%
\textbf{TOR}\xspace
}
\newcommand\Spec{%
\text{Spec}
}

\newcommand\Aut{%
\text{Aut}\hspace{0.05 cm}
}
\newcommand\Mor{%
\text{Mor}%\hspace{0.05 cm}
}
\newcommand\Hom{%
\text{Hom}%\hspace{0.05 cm}
}
\newcommand\sA{%
\mathcal{A}
}

\newcommand\sO{%
\mathcal{O}
}

\newcommand\CoKer{%
\text{CoKer}\hspace{0.05cm}
}

\newcommand\ab{%
\text{ab}
}

\newcommand\sep{%
\text{sep}
}

\newcommand\Ker{%
\text{Ker}\hspace{0.05cm}
}
\newcommand\Gal{%
\text{Gal}
}
\newcommand\rec{%
\text{rec}
}

\newcommand\Res{%
\text{Res}
}

\newcommand\Ind{%
\text{Ind}
}
\newcommand\Inv{%
\text{Inv}
}

\newcommand\id{%
\text{id}
}

\newcommand\tL{%
\text{L}
}
\newcommand\tN{%
\text{N}
}

\newcommand\GL{%
\text{GL}
}
\newcommand\tS{%
\text{S}
}
\newcommand\SL{%
\text{SL}
}
\newcommand\SO{%
\text{SO}
}
\newcommand\sw{%
\text{sw}
}

\newcommand\ra{%
\rightarrow
}
\newcommand\lra{%
\longrightarrow
}

\newcommand\ds{%
\displaystyle
}

\newcommand\otx{%
\hspace{0.03cm}\otimes\hspace{0.03cm}
}

\newcommand\ots[1]{%
\hspace{0.03cm}\otimes_{#1}\hspace{0.03cm}
}
\newcommand\charc{%
\text{char}%\hspace{0.05cm}
}

\newcommand\OP{%
\text{OP}
}
\newcommand\vol{%
\text{vol}%\hspace{0.05cm}
}

\newcommand\un[1]{%
\underline{#1}\xspace
}
\newcommand\ov[1]{%
\overline{#1}
}
\newcommand\acdot{%
\abs{\hspace{0.05cm} \cdot \hspace{0.05cm}}
}

\newcommand\pp{%
^{\prime\prime}
}

\newcommand\hsp{%
\hspace{0.05cm}
}

\usepackage{scalerel,stackengine}
\stackMath
\newcommand\reallywidehat[1]{%
\savestack{\tmpbox}{\stretchto{%
  \scaleto{%
    \scalerel*[\widthof{\ensuremath{#1}}]{\kern-.6pt\bigwedge\kern-.6pt}%
    {\rule[-\textheight/2]{1ex}{\textheight}}%WIDTH-LIMITED BIG WEDGE
  }{\textheight}% 
}{0.5ex}}%
\stackon[1pt]{#1}{\tmpbox}%
}
\DeclareFontFamily{U}{wncy}{}
    \DeclareFontShape{U}{wncy}{m}{n}{<->wncyr10}{}
    \DeclareSymbolFont{mcy}{U}{wncy}{m}{n}
    \DeclareMathSymbol{\Sh}{\mathord}{mcy}{"58} 
\makeatletter
\DeclareRobustCommand\widecheck[1]{{\mathpalette\@widecheck{#1}}}
\def\@widecheck#1#2{%
    \setbox\z@\hbox{\m@th$#1#2$}%
    \setbox\tw@\hbox{\m@th$#1%
       \widehat{%
          \vrule\@width\z@\@height\ht\z@
          \vrule\@height\z@\@width\wd\z@}$}%
    \dp\tw@-\ht\z@
    \@tempdima\ht\z@ \advance\@tempdima2\ht\tw@ \divide\@tempdima\thr@@
    \setbox\tw@\hbox{%
       \raise\@tempdima\hbox{\scalebox{1}[-1]{\lower\@tempdima\box
\tw@}}}%
    {\ooalign{\box\tw@ \cr \box\z@}}}
\makeatother

\newtheoremstyle{xx}% name of the style to be used
  {4pt}% measure of space to leave above the theorem. E.g.: 3pt
  {0pt}% measure of space to leave below the theorem. E.g.: 3pt
  {\upshape}% name of font to use in the body of the theorem
  {\bfseries}% name of head font
  {}% punctuation between head and body
  { }% space after theorem head; " " = normal interword space \footnotesize
  {}
  
\makeatletter 
 \newtheoremstyle{myu}%
  {\upshape\item[ \indent\indent\bf\underline{\theorem@headerfont ##2:}]}%
\makeatother
\makeatletter 
 \newtheoremstyle{myn}%
  {\item[\hskip\labelsep \ \bf ##1 \theorem@headerfont ##2.]}%
\makeatother
\theoremstyle{myn}
\newtheorem{theoremn}{Theorem} %<-- Normal Theorem Definition
\theoremstyle{myu}
{\upshape}
\newtheorem{x}[theoremn]{}%<-- Underlined Theorem Definition

\title{\textbf{Abelian Theory}}
\author{Garth Warner\\
Department of Mathematics\\
University of Washington}
\date{}									% Activate to display a given date or no date
\titleformat{\chapter}[display]
{\normalfont\filcenter\huge\bfseries}{}{0pt}{\large}
\titleformat{\chapter}[display]
{\normalfont\filcenter\huge\bfseries}{}{0pt}{\large}
\usepackage[OT2,OT1]{fontenc} 
\newcommand\cyr
{
\renewcommand\rmdefault{wncyr} 
\renewcommand\sfdefault{wncyss} 
\renewcommand\encodingdefault{OT2} 
\normalfont
\selectfont
}

\DeclareTextFontCommand{\textcyr}{\cyr}

\makeindex
\begin{document}

\maketitle                              % Print title page.

\titlespacing*{\chapter}{0pt}{-50pt}{40pt}
\setlength{\parskip}{0.1em}
\include{_tocX}

\include{_Preface}

\pagenumbering{bychapter}
%\chapter{}
%\newpage
%\setcounter{page}{1}
%\renewcommand{\thepage}{1-\arabic{page}}
%\renewcommand{\thepage}{1-\arabic{page}}
\setcounter{chapter}{0}
\include{_01_GROUP_SCHEMES}

\include{_02_SCH_K}

\include{_03_AFFINE_GROUP_k_SCHEMES}

\include{_04_ALGEBRAIC_TORI}

\include{_05_THE_LLC}

\include{_06_TAMAGAWA_MEASURES}

%\newpage
\include{_references}

\setcounter{page}{1}
\renewcommand{\thepage}{Index-\arabic{page}}
\backmatter
\bibliography{}
\printindex
\end{document}

%% file: _tocX.tex
\begingroup%%----------------------------------->>
\begingroup
\center {\textbf{CONTENTS}}\\
\endgroup
\vspace{1.cm}

%%\qquad\qquad CONTENTS\\
 %%dmc00A_\vspace{2.cm}

\qquad $\S1.\ $ \qquad GROUP SCHEMES\\
%%dmc00A_\vspace{0.3cm}

\qquad $\S2.\ $ \qquad \un{SCH}/$\bk$\\
%%dmc00A_\vspace{0.3cm}

\qquad $\S3.\ $ \qquad AFFINE GROUP $\bk$-SCHEMES\\
%%dmc00A_\vspace{0.3cm}

\qquad $\S4.\ $ \qquad ALGEBRAIC TORI\\
%%dmc00A_\vspace{0.3cm}

\qquad $\S5.\ $ \qquad THE LOCAL LANGLANDS CORRESPONDENCE (LLC)\\
%%dmc00A_\vspace{0.3cm}

\qquad $\S6.\ $ \qquad TAMAGAWA MEASURES\\
%%dmc00A_\vspace{0.3cm}

%%dmc00A_\vspace{0.3cm}
\endgroup %%------------------------------------<<

%% file: _Preface.tex
\begingroup
\center {\textbf{INTRODUCTION}}\\
\endgroup
\vspace{0.5cm}

\ \indent 
This exposition begins with a systematic account of the theory of group schemes, ultimately specializing to algebraic tori.  
This done, I then take up two central objectives.
\vspace{0.3cm}

\qquad\textbullet \ \ 
First, I define the local Langlands correspondece in the case of a $\K$-torus ($\K$ a non-archimedean local field) and prove it when 
\[
T 
\ = \ 
\Res_{\LL / \K} (G_{m, \LL}),
\]
$\LL$ a finite Galois extension of $\K$ ($\Res_{\LL / \K}$ the functor of restriction of the scalars).  
Historically, this was first done by Langlands in the late 1960's and served to confirm his far reaching conjectures for reductive groups in the ``simplest situation'' ($\GL_1$ ``is'' local class field theory \ldots).
\vspace{0.3cm}

\qquad\textbullet \ \ 
Second, working within the context of a $\Q$-torus \mT, I define precisely the terms 
``Tamagawa measure'' 
and 
``Tamagawa number'',
make some computations to illustrate the general picture, and set up Ono's celebrated theorem connecting these notions with the Tate-Shafarevich group.

%% file: _01_GROUP_SCHEMES.tex
\chapter{
$\boldsymbol{\S}$\textbf{1}.\quad  GROUP SCHEMES}
\setlength\parindent{2em}
\setcounter{theoremn}{0}
%%----------------------------------------------------------------------------------------------01

\ \indent 

\begin{x}{\small\bf NOTATION} \ %01
\bSCH is the category of schemes, \bRNG is the category of commutative rings with unit.
\end{x}
\vspace{0.3cm}

Fix a scheme \mS $-$then the category 
$\bSCH/S$
\index{$\bSCH/S$} 
of 
\un{schemes over \mS}
\index{schemes over \mS} 
(or of \un{\mS-schemes})\index{\mS-schemes} 
is the category whose objects are the morphisms $X \ra S$ of schemes and whose morphisms 
\[
\Mor(X \ra S, Y \ra S)
\]
are the morphisms $X \ra Y$ of schemes with the property that the diagram 
\[
\begin{tikzcd}[sep=large]
X  \ar{d} \ar{rr} &&{Y} \ar{d}\\
S 
\arrow[rr,shift right=0.5,dash] \arrow[rr,shift right=-0.5,dash] 
&&{S}
\end{tikzcd}
\]
commutes.

[Note: \ Take $S =\Spec(\Z)$ $-$then 
\[
\bSCH / S \ =  \ \bSCH.]
\]

\begin{x}{\small\bf \un{N.B.}} \ %02
If $S = \Spec(A)$ (\mA in \bRNG) is an affine scheme, then the terminology is ``schemes over \mA'' (or ``\mA-schemes'') and one writes 
$\bSCH / A$ in place of $\bSCH / \Spec(A)$.
\end{x}
\vspace{0.3cm}

\begin{x}{\small\bf NOTATION} \ %03
Abbreviate $\Mor(X \ra S, Y\ra S)$ to $\Mor_S(X,Y)$ (or to ($\Mor_A(X,Y)$ if $S = \Spec(A)$).
\end{x}
\vspace{0.3cm}

\begin{x}{\small\bf REMARK} \ %04
The \mS-scheme $\id_S:S \ra S$ is a final object in $\bSCH / S$.
\end{x}
\vspace{0.3cm}

%%----------------------------------------------------------------------------------------------02

\begin{x}{\small\bf THEOREM} \ %05
$\bSCH / S$ has pullbacks:
\[
\begin{tikzcd}[sep=large]
{X \times_S Y}  \ar{d} \ar{rr} &&{Y} \ar{d}\\
{X}
\ar{rr}
&&{S}
\end{tikzcd}
.
\]

[Note: \ Every diagram 
\[
\begin{tikzcd}%[sep=large]
{Z}  \ar[bend right]{dddr}[swap]{u}  \arrow[drrr, bend left, "v"] \\
&{X \times_S Y} \ar{dd}[swap]{p} \ar{rr}{q} &&{Y} \ar{dd}{q}\\
\\
&{X} \ar{rr}[swap]{f} &&{S}
\end{tikzcd}
\quad (f \circ u = g \circ v)
\]
admits a unique filler 
\[
(u,v)_S:Z \ra X \times_S Y
\]
such that 
$
\begin{cases}
\ p \circ (u,v)_S \ = \ u\\
\ q \circ (u,v)_S \ = \ v
\end{cases}
.]
$
\end{x}
\vspace{0.3cm}

\begin{x}{\small\bf FORMALITIES} \ %06
Let $X, Y, Z$ be objects in $\bSCH / S$ $-$then 
\begin{align*}
&X \times_S S \approx X,\\
&X \times_S Y \approx Y \times_S X,
\end{align*}
and 
\[
(X \times_S Y) \times_S Z \ \approx \  X \times_S (Y \times_S Z).
\]
\end{x}
\vspace{0.3cm}
%%----------------------------------------------------------------------------------------------03

\begin{x}{\small\bf REMARK} \ %07
If $X, Y, X^\prime, Y^\prime$ are objects in $\bSCH / S$ and if 
$u:X \ra X^\prime$, $v:Y \ra Y^\prime$ are $S$-morphisms, then there is a unique morphism 
$u \times_S v$ (or just $u \times v$) rendering the diagram
\[
\begin{tikzcd}%[sep=large]
{X}   \ar{rr}{u}
&&{X^\prime} \ar{rr} 
&&{S} \arrow[dd,shift right=0.5,dash] \arrow[dd,shift right=-0.5,dash] \\
{X \times_S Y} \ar{d}[swap]{q} \ar{u}{p} \ar{rr}{u \times_S v} 
&&{X^\prime \times_S Y^\prime} \ar{u} \ar{d}\\
{Y} \ar{rr}[swap]{v} 
&&{Y^\prime} \ar{rr} 
&&{S}
\end{tikzcd}
\]
commutative.

[Spelled out, 
\[
u \times_S v \ = \ (u \circ p, v \circ q)_S.]
\]
\end{x}
\vspace{0.3cm}

\begin{x}{\small\bf BASE CHANGE} \ %08
Let $u:S^\prime \ra S$ be a morphism in \bSCH.
\vspace{0.2cm}

\qquad \textbullet \quad If $X \ra S$ is an $S$-object, then $X \times_S S^\prime$ is an $S^\prime$-object via the projection 
\[
X \times_S S^\prime \lra S^\prime, 
\]
denoted $u^*(X)$ or $X_{(S^\prime)}$ and called the 
\un{base change of \mX by $u$}.
\index{base change of \mX by $u$}
\vspace{0.2cm}

\qquad \textbullet \quad  
If $X \ra S$, $Y \ra S$ are $S$-objects and if $f:(X \ra S) \ra (Y \ra S)$ is an $S$-morphism, then 
\[
\begin{tikzcd}[sep=large]
{X \times_S S^\prime}  \ar{d} \ar{rrr}{f \times_S \id_{S^\prime}}
&&&{Y \times_S S^\prime} \ar{d}\\
{S^\prime}
\arrow[rrr,shift right=0.5,dash] \arrow[rrr,shift right=-0.5,dash] 
&&&{S^\prime}
\end{tikzcd}
\]
is a morphism of $S^\prime$-objects, denoted $u^*(f)$ or $f_{(S^\prime)}$ and called 
\un{the base change of $f$ by $u$}.
\index{the base change of $f$ by $u$}
\end{x}
\vspace{0.3cm}
%%----------------------------------------------------------------------------------------------04

These considerations thus lead to a functor 
\[
u^*: \bSCH / S \ra \bSCH /  S^\prime
\]
called \un{the base change by $u$}.
\index{the base change by $u$}
\vspace{0.3cm}

\begin{x}{\small\bf \un{N.B.}} \ %09
If $u^\prime:S\pp \ra S^\prime$ is another morphism in \bSCH, then the functors 
$(u \circ u^\prime)^*$ and $(u^\prime)^* \circ u$ from $\bSCH / S$ to $\bSCH /  S\pp$ are isomoprhic. 
\end{x}
\vspace{0.3cm}

\begin{x}{\small\bf LEMMA} \ %10
Let $u:S^\prime \ra S$ be a morphism in \bSCH.  
Suppose that $T^\prime \ra S^\prime$ is an $S^\prime$-object $-$then $T^\prime$ can be viewed as an $S$-object $T$ 
via postcomposition with $u$ and there are canonical mutually inverse bijections
\[
\begin{tikzcd}[sep=small]
{\Mor_{S^\prime}(T^\prime,X_{(S^\prime)})} \arrow[r,shift left=0.97] 
&{\Mor_S(T,X)} \arrow[l,shift left=0.97]
\end{tikzcd}
\]
functorial in $T^\prime$ and \mX.
\end{x}
\vspace{0.3cm}

\begin{x}{\small\bf NOTATION} \ %11
Each $S$-scheme $X \ra S$ determines a functor 
\[
(\bSCH /  S)^\OP \ra \bSET,
\]
viz. the assignment
\[
T \ra \Mor_S(T,X) \ \equiv \ X_S(T),
\]
the set of \un{\mT-valued points of \mX}.
\index{\mT-valued points of \mX}
\vspace{0.1cm}

[Note: \ In terms of category theory, 
\[
X_S(T) \ = \ h_{X \ra S} (T \ra S).]
\]
\end{x}
\vspace{0.3cm}

%%----------------------------------------------------------------------------------------------05
\begin{x}{\small\bf LEMMA} \ %12
To give a morphism 
$(X \ra S) \overset{f}{\lra} (Y \ra S)$ in $\bSCH / S$ is equivalent to giving for all $S$-schemes \mT a map
\[
f(T):X_S(T) \ra Y_S(T)
\]
which is functorial in \mT, i.e., for all morphisms 
$u:T^\prime \ra T$ of $S$-schemes the diagram
\[
\begin{tikzcd}[sep=large]
{X_S(T)}  \ar{d}[swap]{X_S(u)} \ar{rr}{f(T)}
&&{Y_S(T)}  \ar{d}{Y_S(u)}\\
{X_S(T^\prime)} \ar{rr}[swap]{f(T^\prime)}
&&{Y_S(T^\prime)}
\end{tikzcd}
\]
commutes.
\end{x}
\vspace{0.3cm}

\begin{x}{\small\bf DEFINITION} \ %13
A \un{group scheme} over \mS (or an \un{\mS-group}) is an object \mG of $\bSCH / S$ and \mS-morphisms
\index{group scheme}
\index{\mS-group}
\begin{align*}
&m:G \times_S G \ra G \qquad \ \  \text{(``multiplication'')}\\
&e:S \ra G \qquad\qquad\qquad \text{(``unit'')}\\ 
&i:G \ra G \qquad\qquad\qquad \text{(``inversion'')}
\end{align*}
such that the diagrams
\[
\begin{tikzcd}[sep=large]
{G \times_S G \times_S G}  \ar{d}[swap]{\id_G \times m} \ar{rr}{m \times \id_G}
&&{G \times_S G}  \ar{d}{m}\\
{G \times_S G} \ar{rr}[swap]{m}
&&{G}
\end{tikzcd}
\]
%%----------------------------------------------------------------------------------------------06
\[
\begin{tikzcd}[sep=large]
{G \times_S S}  \arrow[d,shift right=0.5,dash] \arrow[d,shift right=-0.5,dash] \ar{rr}{(\id_G,e)_S}
&&{G \times_S G}  \ar{d}{m}\\
{G} \ar{rr}[swap]{\id_G}
&&{G}
\end{tikzcd}
\]
\[
\begin{tikzcd}[sep=large]
{G}  \ar{d} \ar{rr}{(\id_G,i)_S}
&&{G \times_S G}  \ar{d}{m}\\
{S} \ar{rr}[swap]{e}
&&{G}
\end{tikzcd}
\]
commute.
\end{x}
\vspace{0.3cm}

\begin{x}{\small\bf REMARK} \ %14
To say that $(G; m,e,i)$ is a group scheme over \mS amounts to saying that \mG is a group object in $\bSCH / S$.
\end{x}
\vspace{0.3cm}

\begin{x}{\small\bf LEMMA} \ %15
Let \mG be an \mS-scheme $-$then \mG gives rise to a group scheme over \mS iff for all \mS-schemes \mT, 
the set $G_S(T)$ carries the structure of a group which is functorial in \mT (i.e., for all \mS-morphisms 
$T^\prime \ra T$, the induced map $G_S(T) \ra G_S(T^\prime)$ is a homomorphism of groups).
\end{x}
\vspace{0.3cm}

\begin{x}{\small\bf REMARK} \ %16
It suffices to define functorial group structures on the $G_S(A)$, where $\Spec(A) \ra S$ is an affine \mS-scheme.

[This is because morphisms of schemes can be ``glued''.]
\end{x}
\vspace{0.3cm}

\begin{x}{\small\bf LEMMA} \ %17
Let $u:S^\prime \ra S$ be a morphism in \bSCH.  Suppose that $(G;m,e,i)$ is 
%%----------------------------------------------------------------------------------------------07
a group scheme over \mS $-$then 
\[
(G \times_S S^\prime; m_{(S^\prime)}, e_{(S^\prime)}, i_{(S^\prime)})
\]
is a group scheme over $S^\prime$.
\vspace{0.2cm}

[Note: \ For every $S^\prime$-object $T^\prime \ra S^\prime$, 
\[
(G \times_S S^\prime)_{S^\prime}(T^\prime) \ = \ G_S(T),
\]
where \mT is the \mS-object $T^\prime \ra S^\prime \overset{u}{\lra} S$.]
\end{x}
\vspace{0.3cm}

\begin{x}{\small\bf THEOREM} \ %18
If $(X,\sO_X)$ is a locally ringed space and if \mA is a commutative ring with unit, then there is a functorial set-theoretic bijection 
\[
\Mor(S,\Spec(A)) \ \approx \ \Mor(A,\Gamma(X,\sO_X)).
\]
\vspace{0.2cm}

[Note: \ 
The ``Mor'' on the LHS is in the category of locally ringed spaces and the ``Mor'' on the RHS is in the category of commutative rings with unit.]
\end{x}
\vspace{0.3cm}

\begin{x}{\small\bf EXAMPLE} \ %19
Take $S = \Spec(\Z)$ and let 
\[
A^n \ = \ \Spec(\Z[t_1, \ldots, t_n]).
\]
Then for every scheme \mX,
\begin{align*}
\Mor(X,A^n) \ 
&\approx\ \Mor(\Z[t_1, \ldots, t_n], \ \Gamma(X,\sO_X))\\
&\approx\ \Gamma(X,\sO_X)^n \qquad (\phi \ra (\phi(t_1), \ldots, \phi(t_n))).
\end{align*}
Therefore $A^n$ is a group object in \bSCH called 
\un{affine $n$-space}.
\index{affine $n$-space}
\vspace{0.2cm}

[Note: \ 
Here $\Gamma(X,\sO_X)$ is being viewed as an additive group, hence the underlying multiplicative structure is being ignored.]
\end{x}
\vspace{0.3cm}

%%----------------------------------------------------------------------------------------------08

\begin{x}{\small\bf \un{N.B.}} \ %20
Given any scheme \mS, 
\[
A_S^n \ = \ A^n \times_\Z S \ra S
\]
is an \mS-scheme and for every morphism $S^\prime \ra S$, 
\[
A_S^n \times_S S^\prime \ \approx \ A^n \times_\Z S \times_S S^\prime \ \approx \ A_{S^\prime}^n.
\]
\end{x}
\vspace{0.1cm}

\begin{x}{\small\bf NOTATION} \ %21
Write $G_a$ in place of $A^1$.
\end{x}
\vspace{0.3cm}

\begin{x}{\small\bf NOTATION} \ %22
Given \mA in \bRNG, denote
\[
G_a \times_\Z \Spec(A)
\]
by $G_a {\otimes} A$ or still, by $G_{a,A}$.
\end{x}
\vspace{0.3cm}

\begin{x}{\small\bf \un{N.B.}} \ %23
\begin{align*}
G_{a,A} \ 
&=\ \Spec(\Z[t]) \times_\Z \Spec(A)\\
&=\ \Spec(\Z[t] {\otimes} A) \\
&=\ \Spec(A[t]).
\end{align*}
\end{x}
\vspace{0.3cm}

\begin{x}{\small\bf LEMMA} \ %24
$G_{a,A}$ is a group object in $\bSCH / A$.
\end{x}
\vspace{0.3cm}

There are two other ``canonical'' examples of group objects in $\bSCH / A$.
\vspace{0.2cm}

\qquad \textbullet \quad $G_{m,A} = \Spec(A[u,v] / (uv - 1))$\\
which assigns to an \mA-scheme \mX the multiplicative group $\Gamma(X,\sO_X)^\times$ 
of invertible elements in the ring $\Gamma(X,\sO_X)$.
\vspace{0.2cm}

\qquad \textbullet \quad $\GL_{n,A} = \Spec(A[t_{11}, \ldots, t_{nn}, \det(t_{ij})^{-1}])$\\
%%----------------------------------------------------------------------------------------------09
which assigns to an \mA-scheme \mX the group 
\[
\GL_n(\Gamma(X,\sO_X))
\]
of invertible $n \times n$-matrices with entries in the ring $\Gamma(X,\sO_X)$.
\vspace{0.3cm}

\begin{x}{\small\bf DEFINITION} \ %25
If \mG and \mH are \mS-groups, then a 
\un{homomorphism from \mG to \mH}
\index{homomorphism of \mS-groups} 
is a morphism 
$f:G \ra H$ of \mS-schemes such that for all \mS-schemes \mT the induced map 
$f(T):G_S(T) \ra H_S(T)$ is a group homomorphism.
\end{x}
\vspace{0.3cm}

\begin{x}{\small\bf EXAMPLE} \ %26
Take $S = \Spec(A)$ $-$then 
\[
{\det}_A : \GL_{n,A} \ra G_{m,A}
\]
is a homomorphism.
\end{x}
\vspace{0.3cm}

\begin{x}{\small\bf DEFINITION} \ %27
Let \mG be a group scheme over \mS $-$then a subscheme (resp. an open subscheme, resp. a closed subscheme) 
$H \subset G$ is called an \mS-subgroup scheme (resp. an open \mS-subgroup scheme, resp. a closed \mS-subgroup scheme) 
if for every \mS-scheme \mT, $H_S(T)$ is a subgroup of $G_S(T)$.
\end{x}
\vspace{0.3cm}

\begin{x}{\small\bf EXAMPLE} \ %28
Given a positive integer $n$, $\un{\mu}_{n,A}$ is the group object in $\bSCH /A$ which assigns to an \mA-scheme \mX 
the multiplicative subgroup of $\Gamma(X,\sO_X)^\times$ consisting of those $\phi$ such that $\phi^n = 1$, thus 
\[
\un{\mu}_{n,A} \ = \ \Spec(A[t]/(t^n - 1))
\]
and $\un{\mu}_{n,A}$  is a closed \mA-subgroup of $G_{m,A}$.
\end{x}
\vspace{0.3cm}

\begin{x}{\small\bf EXAMPLE} \ %29
Fix a prime number $p$ and suppose that \mA has characteristic $p$.\\
\end{x}
%%----------------------------------------------------------------------------------------------10
Given a positive integer $n$, $\un{\alpha}_{n,A}$ is the group object in $\bSCH /  A$ which assigns to an \mA-scheme 
\mX the additive subgroup of $\Gamma(X,\sO_X)$ consisting of those $\phi$ such that $\phi^{p^n} = 0$, thus
\[
\un{\alpha}_{n,A} \ = \ \Spec(A[t] / t^{p^n}))
\]
and $\un{\alpha}_{n,A}$ is a closed \mA-subgroup of $G_{a,A}$.
\vspace{0.3cm}

\begin{x}{\small\bf CONSTRUCTION} \ %30
Let $f:G \ra H$ be a homomorphism of \mS-groups.  Define $\Ker(f)$ by the pullback square
\[
\begin{tikzcd}[sep=large]
{\Ker(f) = S \times_H G}  \ar{d} \ar{rr}
&&{G}  \ar{d}{f}\\
{S} \ar{rr}[swap]{e}
&&{H}
\end{tikzcd}
.
\]
Then for all \mS-schemes \mT, 
\[
\begin{tikzcd}[sep=large]
{\Mor_S(T,\Ker(f)) \ = \ \Ker(G_S(T)} \ar{r}{f(T)} &{H_S(T))}
\end{tikzcd}
,
\]
so $\Ker(f)$ is an \mS-group.
\end{x}
\vspace{0.1cm}

\begin{x}{\small\bf EXAMPLE} \ %31
The kernel of $\det_A$ is $\SL_{n,A}$.
\end{x}
\vspace{0.1cm}

\begin{x}{\small\bf \un{N.B.}} \ %32
Other kernels are $\un{\mu}_{n,A}$ and $\un{\alpha}_{n,A}$.
\end{x}
\vspace{0.1cm}

\begin{x}{\small\bf CONVENTION} \ %33
If P is a property of morphisms of schemes, then an \mS-group \mG has property P if this is the case of its structural morphism 
$G \ra S$.
\end{x}
\vspace{0.3cm}

%%----------------------------------------------------------------------------------------------11
E.g.: \ The property of morphisms of schemes being quasi-compact, locally of finite type, separated, \'etale etc.

\begin{x}{\small\bf LEMMA} \ %34
Let 
\[
\begin{tikzcd}[sep=large]
{X^\prime}  \ar{d}[swap]{f^\prime} \ar{rr}
&&{X}  \ar{d}{f}\\
{Y^\prime} \ar{rr}
&&{Y}
\end{tikzcd}
\]
be a pullback square in \bSCH.  
Suppose that $f$ is a closed immersion $-$then the same holds for $f^\prime$.
\end{x}
\vspace{0.3cm}

\begin{x}{\small\bf APPLICATION} \ %35
Let $g:Y \ra X$ be a morphism of schemes that has a section $s:X \ra Y$.  
Assume: $g$ is separated $-$then $s$ is a closed immersion.

[The commutative diagram
\[
\begin{tikzcd}[sep=large]
{X}  \ar{d}[swap]{s} \ar{rr}{s}
&&{Y}  \ar{d}{\Delta_{Y/X}}\\
{Y} \ar{rr}[swap]{(\id_Y,s \circ g)_X}
&&{Y \times_X Y}
\end{tikzcd}
\]
is a pullback square in \bSCH.  But $g$ is separated, hence the diagonal morphism 
$\Delta_{Y/X}$ is a closed immersion.  Now quote the preceding lemma.]
\end{x}
\vspace{0.3cm}

If $G \ra S$ is a group scheme over \mS, then the composition 
\[
S \overset{e}{\lra} G \lra S
\]
is $\id_S$.  Proof: \ $e$ is an \mS-morphism and the diagram
%%----------------------------------------------------------------------------------------------12
\[
\begin{tikzcd}[sep=large]
{S}  \ar{d}[swap]{\id_s} \ar{rr}{e}
&&{G}  \ar{d}\\
{S} \arrow[rr,shift right=0.5,dash] \arrow[rr,shift right=-0.5,dash] 
&&{S}
\end{tikzcd}
\]
commutes.  Therefore $e$ is a section for the structural morphism $G \ra S$:
\[
G \lra S \overset{e}{\lra} G.
\]
\vspace{0.3cm}

\begin{x}{\small\bf LEMMA} \ %36
Let $G \ra S$ be a group scheme over \mS $-$then the structural morphism $G \ra S$ is separated iff 
$e:S \ra G$ is a closed immersion.
\vspace{0.2cm}

[To see that ``closed immersion'' $\implies$ ``separated'', consider the pullback square
\[
\begin{tikzcd}[sep=large]
{G}  \ar{d}[swap]{\Delta_{G/S}} \ar{rr}
&&{S}  \ar{d}{e}\\
{G \times_S G} \ar{rr}[swap]{m \circ (\id_G \times i)}
&&{G}
\end{tikzcd}
.]
\]
\end{x}
\vspace{0.1cm}

\begin{x}{\small\bf LEMMA} \ %37
If \mS is a discrete scheme, then every \mS-group is separated.
\end{x}
\vspace{0.1cm}

\begin{x}{\small\bf APPLICATION} \ %38
Take $S = \Spec(\bk)$, where \bk is a field $-$then the structural morphism 
$X \ra \Spec(\bk)$ of a \bk-scheme \mX is separated.
\end{x}
\vspace{0.1cm}
%%%%%%%%%%%%%%%%%%%%%%%%%%%%%%%%%%%%%%
%%%%%%%%%%%%%%%%%%%%%%%%%%%%%%%%%%%%%%
%%%%%%%%%%%%%%%%%%%%%%%%%%%%%%%%%%%%%%

%% file: _02_SCH_K.tex
\chapter{
$\boldsymbol{\S}$\textbf{2}.\quad  \un{SCH} / $\bk$}
\setlength\parindent{2em}
\setcounter{theoremn}{0}
%%----------------------------------------------------------------------------------------------01

\ \indent 

Fix a field \bk.

\vspace{0.2cm}

\begin{x}{\small\bf DEFINITION} \ %01
A 
\un{\bk-algebra} 
\index{\bk-algebra}
is an object in \bRNG and a ring homomorphism $\bk \ra A$.
\end{x}
\vspace{0.3cm}

\begin{x}{\small\bf NOTATION} \ %02
$\bALG / \bk$ is the category whose objects are the \bk-algebras $\bk \ra A$ and whose morphisms 
\[
(\bk \ra A) \lra (\bk \ra B)
\]
are the ring homomorphisms $A \ra B$ with the property that the diagram 
\[
\begin{tikzcd}[sep=large]
{A}   \ar{rr}
&&{B}\\
{\bk}
\ar{u}
\arrow[rr,shift right=0.5,dash] \arrow[rr,shift right=-0.5,dash] 
&&{\bk} 
\ar{u} 
\end{tikzcd}
\]
commutes.
\end{x}
\vspace{0.3cm}

\begin{x}{\small\bf DEFINITION} \ %03
Let \mA be a \bk-algebra $-$then \mA is 
\un{finitely generated}
\index{$\bk$-algebra \\ finitely generated} 
if there exists a surjective homomorphism 
$\bk[t_1, \ldots, t_n] \ra A$ of \bk-algebras.
\end{x}
\vspace{0.3cm}

\begin{x}{\small\bf DEFINITION} \ %04
Let \mA be a \bk-algebra $-$then \mA is 
\un{finite} 
\index{$\bk$-algebra \\ finite} 
if there exists a surjective homomorphism 
$\bk^n \ra A$ of \bk-modules.
\end{x}
\vspace{0.3cm}

\begin{x}{\small\bf \un{N.B.}} \ %05
A finite \bk-algebra is finitely generated.
\end{x}
\vspace{0.3cm}

Recall now that $\bSCH /  \bk$ stands for $\bSCH / \Spec(\bk)$.
\vspace{0.2cm}

\begin{x}{\small\bf LEMMA} \ %06
The functor 
\[
A \ra \Spec(A)
\]
%%----------------------------------------------------------------------------------------------02
from $(\bALG / \bk)^\OP$ to $\bSCH /  \bk$ is fully faithful.
\end{x}
\vspace{0.3cm}

\begin{x}{\small\bf DEFINITION} \ %07
Let $X \ra \Spec(\bk)$ be a \bk-scheme $-$then \mX is 
\un{locally of finite type} 
\index{locally of finite type} 
if there exists an affine open covering 
$X = \ds\bigcup\limits_{i \in I} U_i$ 
such that for all $i$, 
$U_i = \Spec(A_i)$, where $A_i$ is a finitely generated \bk-algebra.
\end{x}
\vspace{0.3cm}

\begin{x}{\small\bf DEFINITION} \ %08
Let $X \ra \Spec(\bk)$ be a \bk-scheme $-$then \mX is of 
\un{finite type} 
\index{finite type} 
if \mX is locally of finite type and quasi-compact.
\end{x}
\vspace{0.3cm}

\begin{x}{\small\bf LEMMA} \ %09
If a \bk-scheme $X \ra \Spec(\bk)$ is locally of finite type and if $U \subset X$ is an open affine subset, then 
$\Gamma(U, \sO_X)$ is a finitely generated \bk-algebra.
\end{x}
\vspace{0.3cm}

\begin{x}{\small\bf APPLICATION} \ %10 
If \mA is a finitely generated \bk-algebra, then the \bk-scheme $\Spec(A) \ra \Spec(\bk)$ is of finite type.
\end{x}
\vspace{0.3cm}

\begin{x}{\small\bf LEMMA} \ %11
If $X \ra \Spec(\bk)$ is a \bk-scheme of finite type, then all subschemes of \mX are of finite type.
\end{x}
\vspace{0.3cm}

\begin{x}{\small\bf RAPPEL} \ %12
Let $(X,\sO_X)$ be a locally ringed space.  
Given $x \in X$, denote the stalk of $\sO_X$ at $x$ by $\sO_{X,x}$ $-$then 
$\sO_{X,x}$ is a local ring.  And:
\vspace{0.3cm}

\qquad \textbullet \quad $\mfm_x$ is the maximal ideal in $\sO_{X,x}$.
\vspace{0.3cm}

\qquad \textbullet \quad $\kappa(x) = \sO_{X,x} / \mfm_x$ is the residue field of $\sO_{X,x}$.
\end{x}
\vspace{0.3cm}

\begin{x}{\small\bf CONSTRUCTION} \ %13
Let $(X,\sO_X)$ be a scheme.  Given $x \in X$, let $U = \Spec(A)$ be an affine open neighborhood of $x$.  
Denote by $\pe$ the prime ideal of \mA corresponding 
%%----------------------------------------------------------------------------------------------03
to $x$, hence $\sO_{X,x} = \sO_{U,x} = A_\pe$ (the localization of \mA at $\pe$) and the canonical homomorphism 
$A \ra A_\pe$ leads to a morphism 
\[
\Spec(\sO_{X,x}) \ = \ \Spec(A_\pe) \ra \Spec(A) \ = \ U \ \subset X
\]
of schemes (which is independent of the choice of \mU).
\end{x}
\vspace{0.3cm}

\begin{x}{\small\bf \un{N.B.}} \ %14
There is an arrow $\sO_{X,x} \ra \kappa(x)$, thus an arrow 
$\Spec(\kappa(x)) \ra \Spec(\sO_{X,x})$, thus an arrow 
\[
i_x:\Spec(\kappa(x)) \ra X
\]
whose image is $x$.
\end{x}
\vspace{0.3cm}

Let $\K$ be any field, let $f:\Spec(\K) \ra X$ be a morphism of schemes, and let $x$ be the image of the unique point $p$ of 
$\Spec(\K)$.  Since $f$ is a morphism of locally ringed spaces, at the stalk level there is a homomorphism 
\[
\sO_{X,x} \ra \sO_{\Spec(\K),p} \ = \ \K
\]
of local rings meaning that the image of the maximal ideal $\mfm_x \subset \sO_{X,x}$ is contained in the maximal ideal 
$\{1\}$ of $\K$, so there is an induced homomorphism 
\[
\iota: \kappa(x) \ra \K.
\]
Consequently, 
\[
f \ = \ i_x \circ \Spec(\iota).
\]
\vspace{0.2cm}

\begin{x}{\small\bf SCHOLIUM} \ %15
There is a bijection 
\[
\Mor(\Spec(\K),X) \ra \{(x,\iota): x \in X, \ \iota :\kappa(x) \ra \K \}.
\]

\end{x}
\vspace{0.3cm}

%%----------------------------------------------------------------------------------------------04

If $X \ra \Spec(\bk)$ is a \bk-scheme, then for any $x \in X$, there is an arrow 
\[
\Spec(\kappa(x)) \ra X,
\]
from which an arrow 
\[
\Spec(\kappa(x)) \ra \Spec(\bk), 
\]
or still, an arrow $\bk \ra \kappa(x)$.

\vspace{0.2cm}

\begin{x}{\small\bf LEMMA} \ %16
Let $X \ra \Spec(\bk)$ be a \bk-scheme locally of finite type $-$then 
$x \in X$ is closed iff the field extension $\kappa(x)/\bk$ is finite.
\end{x}
\vspace{0.3cm}

\begin{x}{\small\bf APPLICATION} \ %17
Let $X \ra \Spec(\bk)$ be a \bk-scheme locally of finite type.  
Assume: \ \bk is algebraically closed $-$then 
\begin{align*}
\{x \in X: x \ \text{closed}\} \ 
&=\ \{x \in X: \bk = \kappa)x)\}\\
&=\ \Mor_\bk(\Spec(\bk),X) \equiv X(\bk).
\end{align*}
 
\end{x}
\vspace{0.3cm}

\begin{x}{\small\bf DEFINITION} \ %18
A subset \mY of a topological space \mX is 
\un{dense} 
\index{dense} 
in \mX if $\ov{Y} = X$.
\end{x}
\vspace{0.3cm}

\begin{x}{\small\bf DEFINITION} \ %19
A subset \mY of a topological space \mX is 
\un{very dense} 
\index{very dense} 
in \mX if for every closed subset 
$F \subset X$, $\ov{F \hsp \cap \hsp Y} = F$.
\end{x}
\vspace{0.3cm}

\begin{x}{\small\bf \un{N.B.}} \ %20
If \mY is very dense in \mX, then \mY is dense in \mX.  
\vspace{0.2cm}

[Take $F = X: \ov{X \hsp \cap \hsp Y} = \ov{Y} = X$.]
\end{x}
\vspace{0.3cm}

\begin{x}{\small\bf LEMMA} \ %21
Let $X \ra \Spec(\bk)$ be a \bk-scheme locally of finite type $-$then 
\[
\{x \in X: x \ \text{closed}\}
\]
is very dense in \mX.
\end{x}
\vspace{0.3cm}

\begin{x}{\small\bf DEFINITION} \ %22
Let $X \ra \Spec(\bk)$ be a \bk-scheme $-$then a point $x \in X$ is 
%%----------------------------------------------------------------------------------------------05
called 
\un{\bk-rational}
\index{\bk-rational} 
if the arrow $\bk \ra \kappa(x)$ is an isomorphism.
\end{x}
\vspace{0.3cm}

\begin{x}{\small\bf \un{N.B.}} \ %23
Sending a \bk-morphism $\Spec(\bk) \ra X$ to its image sets up a bijection between the set 
\[
X(\bk) \ = \ \Mor_\bk(\Spec(\bk),X)
\]
and the set of \bk-rational points of \mX.
\end{x}
\vspace{0.3cm}

\begin{x}{\small\bf REMARK} \ %24
$X(\bk)$ may very well be empty.
\vspace{0.2cm}

[Consider what happens if $\bk^\prime/\bk$ is a proper field extension.]
\end{x}
\vspace{0.3cm}

Given a \bk-scheme $X \ra \Spec(\bk)$ and a field extension $\K/\bk$, let 
\[
X(\K) \ = \ \Mor_\bk(\Spec(\K),X)
\]
be the set of $\K$-valued points of \mX.  
If $x:\Spec(\K) \ra X$ is a $\K$-valued point with image $x \in X$, then there are field extensions
\[
\bk \ra \kappa(x) \ra \K.
\]

\vspace{0.2cm}

\begin{x}{\small\bf \un{N.B.}} \ %25
$\Spec(\K)$ is a \bk-scheme, the structural morphism 
$\Spec(\K) \ra \Spec(\bk)$ 
being derived from the arrow of inclusion 
$j:\bk \ra \K$.]
\end{x}
\vspace{0.3cm}

Let $G = \Gal(\K/\bk)$.  Given $\sigma:\K \ra \K$ in \mG, 
\[
\Spec(\sigma): \Spec(\K) \ra \Spec(\K),
\]
hence 
\[
\begin{tikzcd}[sep=small]
{\Spec(\K)} \ar{rr}{\Spec(\sigma)} 
&&{\Spec(\K)} \ar{r}{x} 
&{X}
\end{tikzcd}
,
\]
and we put
\[
\sigma \cdot x \ = \ x \hsp \circ \hsp \Spec(\sigma).
\]
%%----------------------------------------------------------------------------------------------06

\qquad \textbullet \quad $\sigma \cdot x$ is a $\K$-valued point.
\vspace{0.3cm}

[There is a commutative diagram
\[
\begin{tikzcd}[sep=large]
{\K}  \ar{rr}{\sigma}
&&{\K}\\
{\bk} \ar{u}{j} \ar{rr}[swap]{\id_\bk}
&&{\bk} \ar{u}[swap]{j}
\end{tikzcd}
,
\]
so $\sigma \circ j = j \circ \id_\bk = j$, and if $\pi:X \ra \Spec(\bk)$ is the structural morphism, there is a commutative diagram 
\[
\begin{tikzcd}[sep=large]
{\Spec(\K)} \ar{rr}{x} \ar{d}[swap]{\Spec(j)}
&&{X} \ar{d}{\pi}\\
{\Spec(\bk)} \arrow[rr,shift right=0.5,dash] \arrow[rr,shift right=-0.5,dash] 
&&{\Spec(\bk)}
\end{tikzcd}
,
\]
so $\pi \circ x = \Spec(j)$.  The claim then is that the diagram 
\[
\begin{tikzcd}[sep=large]
{\Spec(\K)} \ar{d}[swap]{\Spec(j)} \ar{rr}{x \hsp \circ \hsp \Spec(\sigma)} 
&&{X} \ar{d}{\pi}\\
{\Spec(\bk)} \arrow[rr,shift right=0.5,dash] \arrow[rr,shift right=-0.5,dash] 
&&{\Spec(\bk)}
\end{tikzcd}
\]
commutes.  But
\begin{align*}
\pi \circ x \circ \Spec(\sigma) \ 
&=\ \Spec(j) \circ \Spec(\sigma)\\
&=\ \Spec(\sigma \circ j)\\
&=\  \Spec(j).]
\end{align*}

\qquad \textbullet \quad The operation 
\[
\begin{cases}
\ G \times X(\K) \ra X(\K)\\
\ (\sigma,x) \ra \sigma \cdot x
\end{cases}
\]
%%----------------------------------------------------------------------------------------------07
is a left action of \mG on $X(\K)$.
\vspace{0.2cm}

[Given $\sigma, \tau \in G: \K \overset{\tau}{\lra} \K \overset{\sigma}{\lra} \K$, 
it is a question of checking that 
\[
(\sigma \circ \tau) \cdot x \ = \  \sigma \cdot (\tau \cdot x).
\]
But the LHS equals
\[
x \circ \Spec(\sigma \circ \tau) \ = \  x \circ \Spec(\tau) \circ \Spec(\sigma)
\]
while the RHS equals
\[
\tau \cdot x \circ \Spec(\sigma) \ = \ x \circ \Spec(\tau) \circ \Spec(\sigma).]
\]

\vspace{0.3cm}

\begin{x}{\small\bf NOTATION} \ %26
Let 
\[
\K^G \ = \ \Inv(G)
\]
be the invariant field associated with \mG.
\end{x}
\vspace{0.3cm}

\begin{x}{\small\bf LEMMA} \ %27
The set $X(\K)^G$ of fixed points in $X(\K)$ for the left action of \mG on $X(\K)$ coincides with the set $X(\K^G)$.
\end{x}
\vspace{0.3cm}

\begin{x}{\small\bf APPLICATION} \ %28
If $\K$ is a Galois extension of \bk, then 
\[
X(\K)^G \ = \ X(\bk).
\]
Take $\K = \bk^\sep$, thus now $G = \Gal(\bk^\sep/\bk)$.
\end{x}
\vspace{0.3cm}

\begin{x}{\small\bf DEFINITION} \ %29
Suppose given a left action $G \times S \ra S$ of \mG on a set \mS $-$then \mS is called a 
\un{\mG-set} 
\index{\mG-set}
if $\forall \ s \in S$, the \mG-orbit $G \cdot s$ is finite or, equivalently, the stabilizer 
$G_S \subset G$ is an open subgroup of \mG.
\end{x}
\vspace{0.3cm}

\begin{x}{\small\bf EXAMPLE} \ %30
Let $X \ra \Spec(\bk)$ be a \bk-scheme locally of finite type $-$then 
%%----------------------------------------------------------------------------------------------08
$\forall \ x \in X(\bk^\sep)$, the \mG-orbit $G \cdot x$ of $x$ in $X(\bk^\sep)$ is finite, hence $X(\bk^\sep)$ is a \mG-set.
\end{x}
\vspace{0.3cm}

\begin{x}{\small\bf DEFINITION} \ %31
Let $X \ra \Spec(X)$ be a \bk-scheme $-$then \mX is \un{\'etale} if it is of the form 
\[
X \ = \ \coprod\limits_{i \in I} \Spec(\K_i), 
\]
where $I$ is some index set and where $\K_i/\bk$ is a finite separable field extension.
\end{x}
\vspace{0.3cm}

There is a category $\bET/\bk$ whose objects are the \'etale \bk-schemes and there is a category 
$G$-\bSET whose objects are the \mG-sets.  

Define a functor 
\[
\Phi: \bET/\bk \ra \text{$G$-\bSET}
\]
by associating with each \mX in $\bET/\bk$ the set $X(\K^\sep)$ equipped with its left $G$-action.
\vspace{0.3cm}

\begin{x}{\small\bf LEMMA} \ %32
$\Phi$ is an equivalence of categories.
\vspace{0.2cm}

PROOF \ 
To construct a functor 
\[
\Psi : \text{$G$-\bSET} \ra \bET/\bk
\]
such that 
\[
\Psi \circ \Phi \ \approx \  \id_{\bET / \bk} \quad \text{and} \quad \Phi \circ \Psi \ \approx \ \id_{G-\bSET},
\]
take a $G$-set \mS and write it as a union of $G$-orbits, say
\[
S \ \approx \ \coprod\limits_{i \in I} G \cdot s_i.
\]
Let $K_i \supset \bk$ be the finite separable field extension inside $\bk^\sep$ corresponding to 
%%----------------------------------------------------------------------------------------------09
the open subgroup $G_{s_i} \subset G$ and assign to \mS the \'etale \bk-scheme 
$\coprod\limits_{i \in I} \Spec(\K_i)$.  
Proceed \ldots \ .
\end{x}
\vspace{0.3cm}

The foregoing equivalence of categories induces an equivalence between the corresponding categories of group objects: 
\[
\text{\'etale group \bk-schemes $\approx$ $G$-groups},
\]
where a $G$-group is a group which is a $G$-set, the underlying left action being by group automorphisms.
\vspace{0.2cm}

\begin{x}{\small\bf CONSTRUCTION} \ %33
Given a group \mM, let $M_\bk$ be the disjoint union
\[
\coprod\limits_M \Spec(\bk),
\] 
the 
\un{constant group \bk-scheme}, 
\index{constant group \bk-scheme}
thus for any \bk-scheme $X \ra \Spec(\bk)$, 
\[
\Mor_\bk(X,M_\bk)
\]
is the set of locally constant maps $X \ra M$ whose group structure is multiplication of functions.

[The terminology is standard but not the best since if \mM is nontrivial, then 
\[
\Mor_\bk(X,M_\bk) \ \approx M
\]
only if \mX is connected.]
\end{x}
\vspace{0.3cm}

\begin{x}{\small\bf EXAMPLE} \ %34
For any \'etale group \bk-scheme \mX, 
\[
X \times_\bk \Spec(\bk^\sep) \ \approx X(\bk^\sep)_\bk \times_\bk \Spec(\bk^\sep).
\]
\vspace{0.2cm}

[Note: \ Here (and elsewhere), 
\[
\times_\bk \ = \ \times_{\Spec(\bk)}.]
\]
\end{x}
\vspace{0.3cm}

%%----------------------------------------------------------------------------------------------10

\begin{x}{\small\bf RAPPEL} \ %35
An \mA in \bRNG is 
\un{reduced}
\index{reduced} 
if it has no nilpotent elements $\neq 0$ 
(i.e., $\not\exists \ a \neq 0$: $a^n = 0$ $(\exists \ n$)).
\end{x}
\vspace{0.3cm}

\begin{x}{\small\bf DEFINITION} \ %36
A scheme \mX is 
\un{reduced} 
\index{reduced} 
if for any nonempty open subset $U \subset X$, the ring $\Gamma(U,\sO_X)$ is reduced. 
\vspace{0.2cm}

[Note: \ This is equivalent to the demand that all the local rings $\sO_{X,x}$ $(x \in X)$ are reduced.]
\end{x}
\vspace{0.3cm}

\begin{x}{\small\bf DEFINITION} \ %37
Let \mX be a \bk-scheme $-$then \mX is 
\un{geometrically reduced} 
\index{geometrically reduced} 
if for every field extension $\K \supset \bk$, the 
$\K$-scheme $X \times_\bk \Spec(\K)$ is reduced.
\end{x}
\vspace{0.3cm}

\begin{x}{\small\bf LEMMA} \ %38
If \mX is a reduced \bk-scheme, then for every separable field extension $\K/\bk$, the $\K$-scheme 
$X \times_\bk \Spec(\K)$ is reduced.
\end{x}
\vspace{0.3cm}

\begin{x}{\small\bf APPLICATION} \ %39
Assume: \bk is a perfect field $-$then every reduced \bk-scheme \mX is geometrically reduced.
\end{x}
\vspace{0.3cm}

\begin{x}{\small\bf THEOREM} \ %40
Assume: \bk is of characteristic zero.  Suppose that \mX is a group \bk-scheme which is locally of finite type $-$then 
\mX is reduced, hence is geometrically reduced.
\end{x}
\vspace{0.3cm}

%%%%%%%%%%%%%%%%%%%%%%%%%%%%%%%%%%%%%%
%%%%%%%%%%%%%%%%%%%%%%%%%%%%%%%%%%%%%%
%%%%%%%%%%%%%%%%%%%%%%%%%%%%%%%%%%%%%%

%% file: _03_AFFINE_GROUP_k_SCHEMES.tex
\chapter{
$\boldsymbol{\S}$\textbf{3}.\quad  AFFINE GROUP $\bk$-SCHEMES}
\setlength\parindent{2em}
\setcounter{theoremn}{0}
%%----------------------------------------------------------------------------------------------01

\ \indent 

Fix a perfect field \bk.
\vspace{0.2cm}

[Recall that a field \bk is perfect if every field extension of \bk is separable 
(equivalently, $\charc(\bk) = 0$ or $\charc(\bk) = p > 0$ and the arrow $x \ra x^p$ is surjective).]  

\begin{x}{\small\bf DEFINITION} \ %01
An \un{affine group \bk-scheme} is a group \bk-scheme of the form $\Spec(A)$, where \mA is a \bk-algebra.
\end{x}
\vspace{0.3cm}

\begin{x}{\small\bf EXAMPLE} \ %02
\[
G_{a,\bk} \ = \ \Spec(\bk[t])
\]
is an affine group \bk-scheme.
\end{x}
\vspace{0.3cm}

\begin{x}{\small\bf EXAMPLE} \ %03
\[
G_{m,\bk} \ = \ \Spec(\bk[t,t^{-1}])
\]
is an affine group \bk-scheme.
\end{x}
\vspace{0.3cm}

\begin{x}{\small\bf EXAMPLE} \ %04
\[
\un{\mu}_{n,\bk} \ = \ \Spec(\bk[t] / (t^n - 1)) \qquad (n \in \N)
\]
is an affine group \bk-scheme.
\end{x}
\vspace{0.3cm}

There is a category $\bGRP / \bk$ whose objects are the group \bk-schemes and whose morphisms are the morphisms 
$f:X \ra Y$ of \bk-schemes such that for all \bk-schemes \mT the induced map
\[
f(T):\Mor_\bk(T,X) \ra \Mor_\bk(T,Y)
\]
is a group homomorphism.
\vspace{0.3cm}

%%----------------------------------------------------------------------------------------------02
\begin{x}{\small\bf NOTATION} \ %05
\[
\bAFF-\bGRP / \bk
\]
is the full subcategory of $\bGRP / \bk$ whose objects are the affine group \bk-schemes.
\end{x}
\vspace{0.3cm}

\begin{x}{\small\bf NOTATION} \ %06
\[
\bGRP-\bALG / \bk
\]
is the category of group objects in $\bALG / \bk$ and
\[
\bGRP-(\bALG / \bk)^\OP
\]
is the category of group objects in $(\bALG / \bk)^\OP$.
\end{x}
\vspace{0.3cm}

\begin{x}{\small\bf LEMMA} \ %07
The functor 
\[
A \ra \Spec(A)
\]
from $(\bALG / \bk)^\OP$ to $\bSCH / \bk$ is fully faithful and restricts to an equivalence
\[
\bGRP-(\bALG / \bk)^\OP \ra \bAFF-\bGRP / \bk.
\]
\end{x}
\vspace{0.3cm}

\begin{x}{\small\bf REMARK} \ %08
An object in $\bGRP-(\bALG / \bk)^\OP$ 
is a \bk-algebra \mA which carries the structure of a commutative Hopf algebra over \bk: 
$\exists$ \bk-algebra homomorphisms 
\[
\Delta:A \ra A \otimes_\bk A, \quad \varepsilon:A \ra \bk, \quad S:A \ra A
\]
satisfying the ``usual'' conditions.
\end{x}
\vspace{0.3cm}

\begin{x}{\small\bf N.B.} \ %09
There is another way to view matters, viz.  any functor 
$\bALG / \bk \ra \bGRP$ which is representable by a \bk-algebra serves to determine an affine group \bk-scheme 
(and vice versa).  
From this perspective, a morphism 
%%----------------------------------------------------------------------------------------------03
$G \ra H$ of affine group \bk-schemes is a natural transformation of functors, i.e., 
a collection of group homomorphisms $G(A) \ra H(A)$ such that if $A \ra B$ is a \bk-algebra homomorphism, 
then the diagram
\[
\begin{tikzcd}%[sep=large]
{G(A)}   \ar{d} \ar{rr} 
&&{H(A)} \ar{d} \\
{G(B)}   \ar{rr} 
&&{H(B)} 
\end{tikzcd}
\]
commutes.
\vspace{0.2cm}

[Note: \ 
Suppose that 
\[
\begin{cases}
\ G \ = \ h^X \ = \ \Mor(X,-)\\
\ H \ = \ h^Y \ = \ \Mor(Y,-)
\end{cases}
.
\]
Then from Yoneda theory, 
\[
\Mor(G,H) \ \approx \ \Mor(Y,X).]
\]
\end{x}
\vspace{0.3cm}

\begin{x}{\small\bf EXAMPLE} \ %10
$\bk[t,t^{-1}]$ represents $G_{m,\bk}$ and 
\[
\bk[t_{11}, \ldots, t_{nn},\det(t_{ij})^{-1}]
\]
represents $\GL_{n,\bk}$.  
Given any \bk-algebra \mA, the determinant is a group homomorphism
\[
\GL_{n,\bk}(A) \ra G_{m,\bk}(A)
\]
and 
\[
{\det}_\bk \in \Mor(\GL_{n,\bk}, G_{m,\bk}).
\]
\vspace{0.2cm}

[Note: \ 
There is a homomorphism 
\[
\bk[t,t^{-1}] \ra \bk[t_{11}, \ldots, t_{nn},\det(t_{ij})^{-1}]
\]
of \bk-algebras that defines $\det_\bk$.  E.g.: \ If $n = 2$, then the homomorphism in question sends 
$t$ to $t_{11}t_{22} - t_{12}t_{21}$.]
\end{x}
\vspace{0.3cm}

%%----------------------------------------------------------------------------------------------04

\begin{x}{\small\bf PRODUCTS} \ %11
Let
\[
\begin{cases}
\ G = h^X \qquad (X \ \text{in} \ \bALG / \bk)\\
\ H = h^Y \qquad (Y \ \text{in} \ \bALG / \bk)
\end{cases}
\]
be affine group \bk-schemes.  Consider the functor
\[
G \times H: \bALG / \bk \ra \bGRP
\]
defined on objects by
\[
A \ra G(A) \times H(A).
\]
Then this functor is represented by the \bk-algebra $X \ots{\bk} Y$:
\begin{align*}
\Mor(X \ots{\bk} Y, A) \ 
&\approx\ \Mor(X,A) \times \Mor(Y,A)\\
&=\ G(A) \times H(A).
\end{align*}
\end{x}
\vspace{0.3cm}

\begin{x}{\small\bf EXAMPLE} \ %12
Take
\[
\begin{cases}
\ G = G_{m,\R}\\
\ H = G_{m,\R}
\end{cases}
.
\]
Then 
\[
(G_{m,\R} \times G_{m,\R})(\R) \ = \ \R^\times \times \R^\times \ = \ \C^\times
\]
and 
\[
(G_{m,\R} \times G_{m,\R})(\C) \ = \ \C^\times \times \C^\times.
\]
\end{x}
\vspace{0.3cm}

Let $\bk^\prime / \bk$ be a field extension $-$then for any \bk-algebra \mA, the tensor product 
$A \ots{\bk} \bk^\prime$ is a $\bk^\prime$-algebra, hence there is a functor
\[
\bALG / \bk \ra \bALG / \bk^\prime
\]
%%----------------------------------------------------------------------------------------------05
termed \un{extension of the scalars}. On the other hand, every $\bk^\prime$-algebra $B^\prime$ can be regarded 
as a \bk-algebra \mB, from which a functor 
\[
\bALG / \bk^\prime \ra \bALG / \bk
\]
termed \un{restriction of the scalars}.

\vspace{0.2cm}

\begin{x}{\small\bf LEMMA} \ %13
For all \bk-algebras \mA and for all $\bk^\prime$-algebras $B^\prime$, 
\[
\Mor_{\bk^\prime}(A \ots{\bk} \bk^\prime,B^\prime) \  \approx \ \Mor_\bk(A,B).
\]
\end{x}
\vspace{0.3cm}

\begin{x}{\small\bf SCHOLIUM} \ %14
The functor ``extension of the scalars'' is a left adjoint for the functor ``restriction of the scalars''.
\end{x}
\vspace{0.3cm}

Let \mG be an affine group \bk-scheme.  Abusing the notation, denote still by \mG the associated functor 
\[
\bALG / \bk \ra \bGRP.
\]
Then there is a functor 
\[
G_{\bk^\prime}:\bALG / \bk^\prime \ra \bGRP,
\]
namely
\[
G_{\bk^\prime}(A^\prime) \ = \ G(A),
\]
where \mA is $A^\prime$ viewed as a \bk-algebra.
\vspace{0.2cm}

\begin{x}{\small\bf LEMMA} \ %15
$G_{\bk^\prime}$ is an affine group $\bk^\prime$-scheme and the assignment $G \ra G_{\bk^\prime}$ is functorial:

\[
\bAFFGRP / \bk \ra \bAFFGRP / \bk^\prime.
\]
\vspace{0.2cm}

[Note: \ 
Suppose that $G = h^X$ $-$then
\begin{align*}
\Mor_{\bk^\prime} (X \ots{\bk} \bk^\prime, A^\prime) \ 
&\approx \ \Mor_\bk(X,A)\\
&=\ G(A) \\
&=\ G_{\bk^\prime}(A^\prime).
\end{align*}
%%----------------------------------------------------------------------------------------------06 (ish)
Therefore $G_{\bk^\prime}$ is represented by $X \ots{\bk} \bk^\prime$: 
\[
G_{\bk^\prime} \ = \ h^{X \ots{\bk} \bk^\prime}.
\]
Matters can also be interpreted ``on the other side'':
\[
\begin{tikzcd}%[sep=large]
{G_{\bk^\prime} = \Spec(X \times_\bk \bk^\prime) = \Spec(X) \times_\bk \Spec(\bk^\prime)} 
\ar{dd} \ar{rr}
&&{\Spec(\bk^\prime)} \ar{dd}\\
\\
{G = \Spec(X)} 
\ar{rr}
&&{\Spec(\bk)}
\end{tikzcd}
.]
\]

\end{x}
\vspace{0.3cm}

\begin{x}{\small\bf DEFINITION} \ %16
$G_{\bk^\prime}$ is said to have been obtained from \mG by 
\un{extension} \un{of the scalars}.
\index{extension of the scalars}
\end{x}
\vspace{0.3cm}

\begin{x}{\small\bf NOTATION} \ %17
Given an affine group $\bk^\prime$-scheme $G^\prime$, let $G_{\bk^\prime/\bk}$ be the functor 
\[
\bALG / \bk \ \ra \ \bGRP
\]
defined by the rule
\[
A \ \ra \ G^\prime(A \ots{\bk} \bk^\prime).
\]
\vspace{0.2cm}

[Note: \ 
If $\bk^\prime = \bk$, then $G_{\bk^\prime/\bk} = G$.]
\end{x}
\vspace{0.3cm}

\begin{x}{\small\bf THEOREM} \ %18
Assume that $\bk^\prime/\bk$ is a finite field extension $-$then $G_{\bk^\prime/\bk}$ is an affine group 
\bk-scheme and the assignment $G^\prime \ra G_{\bk^\prime/\bk}$ is functorial: 
\[
\bAFFGRP / \bk^\prime \ \ra \ \bAFFGRP / \bk.
\]
\end{x}
\vspace{0.3cm}

\begin{x}{\small\bf DEFINITION} \ %19
$G_{\bk^\prime/\bk}$ is said to have been obtained from $G^\prime$ by 
\un{restriction} \un{of the scalars}.
\index{restriction of the scalars}.
\end{x}
\vspace{0.3cm}

%%----------------------------------------------------------------------------------------------07

\begin{x}{\small\bf LEMMA} \ %20
Assume that $\bk^\prime/\bk$ is a finite field extension $-$then for all affine group \bk-schemes \mH, 
\[
\Mor_\bk(H,G_{\bk^\prime/\bk}) \ \approx \Mor_{\bk^\prime}(H_{\bk^\prime},G^\prime).
\]
\end{x}
\vspace{0.3cm}

\begin{x}{\small\bf SCHOLIUM} \ %21
The functor ``restriction of the scalars'' is a right adjoint for the functor ``extension of the scalars''.

\vspace{0.3cm}

[Accordingly, there are arrows of adjuction
$
\begin{cases}
\ G \ra (G_{\bk^\prime})_{\bk^\prime/\bk}\\
\ (G_{\bk^\prime/\bk})_{\bk^\prime} \ra G^\prime
\end{cases}
.]
$

\end{x}
\vspace{0.3cm}

\begin{x}{\small\bf NOTATION} \ %22
\[
\Res_{\bk^\prime/\bk}:\bAFFGRP/\bk^\prime \ra \bAFFGRP / \bk
\]
is the functor defined by setting
\[
\Res_{\bk^\prime/\bk}(G^\prime) \ = \ G_{\bk^\prime/\bk}.
\]
So, by definition,
\[
\Res_{\bk^\prime/\bk}(G^\prime)(A) \ = \ G^\prime (A \ots{\bk} \bk^\prime).
\]
And in particular:
\[
\Res_{\bk^\prime/\bk}(G^\prime)(\bk) \ = \ G^\prime (\bk \ots{\bk} \bk^\prime) \ = \ G^\prime(\bk^\prime).
\]
\end{x}
\vspace{0.3cm}

\begin{x}{\small\bf EXAMPLE} \ %23
Take $G^\prime = A_{\bk^\prime}^n$ $-$then 
\[
\Res_{\bk^\prime/\bk}(A_{\bk^\prime}^n) \ \approx \ A_k^{nd} \qquad (d = [\bk^\prime:\bk]).
\]
\end{x}
\vspace{0.3cm}

%%----------------------------------------------------------------------------------------------08

\begin{x}{\small\bf EXAMPLE} \ %24
Take $\bk = \R$, $\bk^\prime = \C$, $G^\prime = G_{m,\C}$, and consider
\[
\Res_{\C / \R}  (G_{m,\C}) .
\]
Then 
\[
\Res_{\C / \R}  (G_{m,\C}) (\R) \ = \ \C^\times
\]
and 
\[
\Res_{\C / \R}  (G_{m,\C})(\C) \ = \ \C^\times \times \C^\times.
\]
\vspace{0.2cm}

[Note: \ 
\[
\Res_{\C / \R}  (G_{m,\C})
\]
is not isomorphic to $G_{m,\R}$ (its group of real points is $\R^\times$).]
\end{x}
\vspace{0.3cm}

\begin{x}{\small\bf LEMMA} \ %25
Let $\bk^\prime$ be a finite Galois extension of \bk $-$then 
\[
(\Res_{\bk^\prime/\bk}(G^\prime))_{\bk^\prime} \ \approx \ 
\prod\limits_{\sigma \in \Gal(\bk^\prime/\bk)} \sigma G^\prime.
\]
\vspace{0.2cm}

[Note: \ 
$\forall \ \sigma \in \Gal(\bk^\prime/\bk)$, there is a pullback square
\[
\begin{tikzcd}%[sep=large]
{\sigma G^\prime}   \ar{dd} \ar{rr} 
&&{\Spec(\bk^\prime)} \ar{dd}{\Spec(\sigma)} \\
\\
{G^\prime}   \ar{rr} 
&&{\Spec(\bk^\prime)}
\end{tikzcd}
.]
\]
 
\end{x}
\vspace{0.3cm}

\begin{x}{\small\bf EXAMPLE} \ %26
Take $\bk = \R$, $\bk^\prime = \C$, $G^\prime = G_{m,\C)}$ $-$then 
\begin{align*}
(\Res_{\C / \R} (G_{m,\C}))_\C \ 
&\approx\ G_{m,\C} \times \sigma G_{m,\C}\\
&\approx\ G_{m,\C} \times G_{m,\C}.
\end{align*}
\end{x}
\vspace{0.3cm}

Let \mG be an affine group $\bk$-scheme.
\vspace{0.3cm}
%%----------------------------------------------------------------------------------------------09

\begin{x}{\small\bf DEFINITION} \ %27
A 
\un{character} 
\index{character}
of \mG is an element of 
\[
X(G) \ = \ \Mor_\bk(G,G_{m,\bk}).
\]

Given $\chi \in X(G)$, for every \bk-algebra \mA, there is a homomorphism 
\[
\chi(A):G(A)  \ra G_{m,\bk}(A) = A^\times. 
\]

Given $\chi_1, \chi_2 \in X(G)$, define
\[
(\chi_1 + \chi_2)(A):G(A) \ra G_{m,\bk}(A) = A^\times
\]
by the stipulation 
\[
(\chi_1 + \chi_2)(A)(t) \ = \ \chi_1(A)(t) \chi_2(A)(t),
\]
from which a character $\chi_1 + \chi_2$ of \mG, hence $X(G)$ is an abelian group.
\end{x}
\vspace{0.3cm}

\begin{x}{\small\bf EXAMPLE} \ %28
Take $G = G_{m,\bk}$ $-$then the characters of \mG are the morphisms $G \ra G_{m,\bk}$ of the form 
\[
t \ra t^n \qquad (n \in \Z),
\]
i.e., 
\[
X(G) \ \approx \ \Z.
\]
\end{x}
\vspace{0.3cm}

\begin{x}{\small\bf EXAMPLE} \ %29
Take $G = G_{m,\bk} \times \cdots \times G_{m,\bk}$  ($d$ factors) $-$then the characters of \mG are the 
morphisms $G \ra G_{m,\bk}$ of the form 
\[
(t_1, \ldots, t_d) \ra t_1^{n_1} \cdots t_d^{n_d} \qquad (n_1, \ldots, n_d \in \Z),
\]
i.e., 
\[
X(G) \ \approx \ \Z^d.
\]
\end{x}
\vspace{0.3cm}

%%----------------------------------------------------------------------------------------------10

\begin{x}{\small\bf EXAMPLE} \ %30
Given an abelian group \mM, its group algebra $\bk[M]$ is canonically a \bk-algebra.  
Consider the functor 
$D(M):\bALG/\bk \ra \bGRP$ defined on objects by the rule
\[
A \ra \Mor(M, A^\times).
\]
Then $\forall \ A$, 
\[
\Mor(M,A^\times) \ \approx \ \Mor(\bk[M],A),
\]
so $\bk[M]$ represents $D(M)$ which is therefore an affine group \bk-scheme.  And 
\[
X(D(M)) \ \approx \ M,
\]
the character of $D(M)$ corresponding to $m \in M$ being the assignment 
\begin{align*}
D(M)(A) \ 
&=\ 
\Mor(M, A^\times)
%\\
%&\ \quad 
%f \ra f(m) 
\end{align*}
\begin{tikzcd}[sep=small]
\qquad\qquad\qquad\qquad\qquad\qquad\qquad\qquad\qquad\quad {} 
\ar{rrrr}{f \ra f(m)} &&&&{A^\times = G_{m,\bk}(A).}
\end{tikzcd}
\end{x}
\vspace{0.3cm}

\begin{x}{\small\bf NOTATION} \ %31
Given $\chi^\prime \in X(G^\prime)$, let $\tN_{\bk^\prime/\bk}(\chi^\prime)$ stand for the rule that assigns 
to each \bk-algebra \mA the homomorphism 
\[
G_{\bk^\prime/\bk} \ra G_{m,\bk}(A) \ = \ A^\times
\]
defined by the composition
\vspace{0.2cm}

\qquad\qquad\qquad $G_{\bk^\prime / \bk}(A)  \hspace{.59cm} \lra G^\prime(A \ots{\bk} \bk^\prime)$\\

\qquad\qquad\qquad $G^\prime(A \ots{\bk} \bk^\prime)
\lra \ G_{m,\bk^\prime}(A \ots{\bk} \bk^\prime)^\times = \ (A \ots{\bk} \bk^\prime)^\times$ \quad .\\

\qquad\qquad\qquad $(A \ots{\bk} \bk^\prime)^\times \hspace{.15cm} \lra A^\times$.
\vspace{0.2cm}
\\
Here the first arrow is the canonical isomorphism, the second arrow is $\chi^\prime (A \ots{\bk} \bk^\prime)$, 
and the third arrow is the norm map.
\end{x}
\vspace{0.3cm}

%%----------------------------------------------------------------------------------------------11

\begin{x}{\small\bf LEMMA} \ %32
The arrow 
\[
\chi^\prime \ra \tN_{\bk^\prime/\bk}(\chi^\prime)
\]
is a homomorphism 
\[
X(G^\prime) \ra X(G_{\bk^\prime/\bk})
\]
of abelian groups.
\end{x}
\vspace{0.3cm}

\begin{x}{\small\bf THEOREM} \ %33
The arrow 
\[
\chi^\prime \ra \tN_{\bk^\prime/\bk}(\chi^\prime)
\]
is bijective, hence defines an isomorphism 
\[
X(G^\prime) \ra X(G_{\bk^\prime/\bk})
\]
of abelian groups.
\end{x}
\vspace{0.3cm}

\begin{x}{\small\bf APPLICATION} \ %34
Consider
\[
\Res_{\C/\R}(G_{m,\C}).
\]
Then its character group is isomorphic to the character group of $G_{m,\C}$, i.e., to $\Z$.  Therefore
\[
\Res_{\C/\R}(G_{m,\C})
\]
is not isomorphic to $G_{m,\R} \times G_{m,\R}$.
\end{x}
\vspace{0.3cm}

%%%%%%%%%%%%%%%%%%%%%%%%%%%%%%%%%%%%%%
%%%%%%%%%%%%%%%%%%%%%%%%%%%%%%%%%%%%%%
%%%%%%%%%%%%%%%%%%%%%%%%%%%%%%%%%%%%%%

%% file: _04_ALGEBRAIC_TORI.tex
\chapter{
$\boldsymbol{\S}$\textbf{4}.\quad  ALGEBRAIC TORI}
\setlength\parindent{2em}
\setcounter{theoremn}{0}
%%----------------------------------------------------------------------------------------------01

\ \indent 

Fix a field \bk of characteristic zero.

\begin{x}{\small\bf DEFINITION} \ %01
Let \mG be an affine group \bk-scheme $-$then \mG is \un{algebraic} if its associated representing \bk-algebra \mA 
is finitely generated.
\end{x}
\vspace{0.1cm}

\begin{x}{\small\bf REMARK} \ %02
It can be shown that every algebraic affine group \bk-scheme is isomorphic to a closed subgroup of some $\GL_{n,\bk}$ 
$(\exists \ n)$.
\end{x}
\vspace{0.1cm}

\begin{x}{\small\bf CONVENTION} \ %03
The term \un{algebraic \bk-group} means ``algebraic affine group \bk-scheme''.
\end{x}
\vspace{0.1cm}

\begin{x}{\small\bf \un{N.B.}} \ %04
It is automatic that an algebraic \bk-group is reduced (cf. \S2, \#40), hence is geometrically reduced (cf. \S2, \#39).
\end{x}
\vspace{0.1cm}

\begin{x}{\small\bf LEMMA} \ %05
Assume that $\bk^\prime/\bk$ is a finite field extension $-$then the functor 
\[
\Res_{\bk^\prime/\bk}: \bAFFGRP /\bk^\prime \ra \bAFFGRP /\bk
\]
sends algebraic $\bk^\prime$-groups to algebraic \bk-groups.
\end{x}
\vspace{0.1cm}

Given a finite field extension $\bk^\prime/\bk$, let $\Sigma$ be the set of \bk-embeddings of $\bk^\prime$ into 
$\bk^\sep$ and identify 
$\bk^\prime \ots{\bk} \bk^\sep$ 
with $(\bk^\sep)^\Sigma$ via the bijection 
which takes $x \otx y$ to the string $(\sigma(x) y)_{\sigma \in \Sigma}$.
\vspace{0.2cm}

\begin{x}{\small\bf LEMMA} \ %06
Let $G^\prime$ be an algebraic $\bk^\prime$-group $-$then 
\[
(G_{\bk^\prime / \bk}) \times_\bk \Spec(\bk^\sep) \ \approx \ \prod\limits_{\sigma \in \Sigma} \sigma G^\prime, 
\]
%%----------------------------------------------------------------------------------------------02
where $\sigma G^\prime$ is the algebraic $\bk^\sep$-group defined by the pullback square
\[
\begin{tikzcd}%[sep=large]
{\sigma G^\prime}   \ar{dd} \ar{rr} 
&&{\Spec(\bk^\sep)} \ar{dd}{\Spec(\sigma)} \\
\\
{G^\prime}   \ar{rr} 
&&{\Spec(\bk^\prime)} 
\end{tikzcd}
.
\]

[Note: \ 
To review, the LHS is 
\[
(\Res_{\bk^\prime/\bk)}(G^\prime))_{\bk^\sep}
\]
and the Galois group $\Gal(\bk^\sep/\bk)$ operates on it through the second factor.  
On the other hand, to each pair $(\tau,\sigma) \in \Gal(\bk^\sep/\bk) \times \Sigma$, there corresponds a bijection 
$\sigma G^\prime \ra (\tau \circ \sigma) G^\prime$ 
leading thereby to an action of $\Gal(\bk^\sep/\bk)$ on 
\[
\prod\limits_{\sigma \in \Sigma} \sigma G^\prime.
\]
The point then is that the identification 
\[
(\Res_{\bk^\prime/\bk)}(G^\prime))_{\bk^\sep} \ \approx \ \prod\limits_{\sigma \in \Sigma} \sigma G^\prime
\]
respects the actions, i.e., is $\Gal(\bk^\sep /\bk)$-equivariant.]
\end{x}
\vspace{0.1cm}

\begin{x}{\small\bf \un{N.B.}} \ %07
Consider the commutative diagram 
\[
\begin{tikzcd}%[sep=large]
{(\tau \circ \sigma) G^\prime}   \ar{dd} \ar{rr} 
&&{\Spec(\bk^\sep)} \ar{dd}{\Spec(\tau)} \\
\\
{\sigma G^\prime}  \ar{dd}\ar{rr} 
&&{\Spec(\bk^\sep)} \ar{dd}{\Spec(\sigma)} \\
\\
{G^\prime} \ar{rr} 
&&{\Spec(\bk^\prime)}
\end{tikzcd}
.
\]
%%----------------------------------------------------------------------------------------------03
Then the ``big'' square is a pullback.  
Since this is also the case of the ``small'' bottom square, it follows that the ``small'' upper square is a pullback.
\end{x}
\vspace{0.1cm}

\begin{x}{\small\bf DEFINITION} \ %08
A \un{split \bk-torus} is an algebraic \bk-group \mT which is isomorphic to a finite product of copies of $G_{m,\bk}$.
\end{x}
\vspace{0.1cm}

\begin{x}{\small\bf EXAMPLE} \ %09
The algebraic $\R$-group
\[
\Res_{\C/\R} (G_{m,\C})
\]
is not a split $\R$-torus (cf. \S3, \#24 and \#34).
\end{x}
\vspace{0.1cm}

\begin{x}{\small\bf LEMMA} \ %10
If \mT is a split \bk-torus, then $X(T)$ is a finitely generated free abelian group.
\end{x}
\vspace{0.1cm}

\begin{x}{\small\bf THEOREM} \ %11
The functor 
\[
T \ra X(T)
\]
from the category of split \bk-tori to the category of finitely generated free abelian groups is a contravariant equivalence of categories.
\end{x}
\vspace{0.1cm}

\begin{x}{\small\bf \un{N.B.}} \ %12
$\forall$ \bk-algebra \mA, 
\[
T(A) \ \approx \ \Mor(X(T),A^\times).
\]

[Note: \ 
Explicated, 
\[
T \ \approx \ \Spec(\bk[X(T)]) \qquad (\text{cf.} \ \S3, \ \#30).
\]
Therefore
\allowdisplaybreaks
\begin{align*}
T(A) \ 
&\approx\ \Mor(\Spec(A),T)\\
&\approx\ \Mor(\Spec(A),\Spec(\bk[X(T)])\\
&\approx\ \Mor(\bk[X(T)],A)\\
&\approx\ \Mor(X(T),A^\times).]
\end{align*}
\end{x}
\vspace{0.1cm}

%%----------------------------------------------------------------------------------------------04

\begin{x}{\small\bf DEFINITION} \ %13
A \un{\bk-torus} is an algebraic \bk-group \mT such that 
\[
T_{\bk^\sep} \ = \ T \times_\bk \Spec(\bk^\sep)
\]
is a split $\bk^\sep$-torus.
\end{x}
\vspace{0.1cm}

\begin{x}{\small\bf \un{N.B.}} \ %14
A split \bk-torus is a \bk-torus.
\end{x}
\vspace{0.1cm}

\begin{x}{\small\bf EXAMPLE} \ %15
Let $\bk^\prime/\bk$ be a finite field extension and take $G^\prime = G_{m,\bk^\prime}$ $-$then 
the algebraic \bk-group $G_{\bk^\prime/\bk}$ is a \bk-torus (cf. \#6).
\end{x}
\vspace{0.1cm}

\begin{x}{\small\bf DEFINITION} \ %16
Let \mT be a \bk-torus $-$then a \un{splitting field} for \mT is a finite field extension $\K/\bk$ such that 
$T_\K$ is a split $\K$-torus.
\end{x}
\vspace{0.1cm}

\begin{x}{\small\bf THEOREM} \ %17
Every \bk-torus \mT admits a splitting field which is minimal (i.e., contained in any other splitting field) and Galois.
\end{x}
\vspace{0.3cm}

\begin{x}{\small\bf NOTATION} \ %18
Given a \bk-scheme \mX and a Galois extension $\K/\bk$, the Galois group $\Gal(\K/\bk)$ operates on 
\[
X_\K \ = \ X \times_\bk \Spec(\K)
\]
via the second term, hence $\sigma \ra 1 \otx \sigma$.
\vspace{0.2cm}

[Note: \ 
$1 \otx \sigma$ is a \bk-automorphism of $X_\K$.]
\end{x}
\vspace{0.3cm}

\begin{x}{\small\bf NOTATION} \ %19
Given \bk-schemes \mX, \mY and a Galois extension $\K/\bk$, the Galois group $\Gal(\K/\bk)$ operates on 
$\Mor_\K(X_\K,Y_\K)$ by the prescription
\[
\sigma f \ = \ (1 \otx \sigma) f(1 \otx \sigma)^{-1}.
\]

[Note: \ 
If $f \in \Mor_\K(X_\K,Y_\K)$, then the condition $\sigma f = f$ for all $\sigma \in \Gal(\K/\bk)$ 
%%----------------------------------------------------------------------------------------------05
is equivalent to the condition that $f$ is the lift of a \bk-automorphism $\phi:X \ra Y$, i.e., $f = \phi \otx 1$.]

\end{x}
\vspace{0.1cm}

\begin{x}{\small\bf LEMMA} \ %20
Let $\K/\bk$ be a Galois extension and let $G = \Gal(\K/\bk)$ $-$then for any \bk-algebra \mA and for any \bk-scheme \mX, 
\[
X(A \ots{\bk} \K)^G \ = \ X(A).
\]

[Note: \ 
This generalizes \S2, \#28 to which it reduces if $A = \bk$.]
\end{x}
\vspace{0.1cm}

\begin{x}{\small\bf DEFINITION} \ %21
Let \mG be a finite group $-$then a \un{\mG-module} is an ablelian group \mM supplied with a homomorphism 
$G \ra \Aut(M)$.
\end{x}
\vspace{0.1cm}

\begin{x}{\small\bf \un{N.B.}} \ %22
A \mG-module is the same thing as a $\Z[G]$-module (in the usual sense when $\Z[G]$ is viewed as a ring).
\end{x}
\vspace{0.1cm}

\begin{x}{\small\bf DEFINITION} \ %23
Let \mG be a finite group $-$then a \un{\mG-lattice} is a $\Z$-free \mG-module \mM of finite rank.
\end{x}
\vspace{0.1cm}

\begin{x}{\small\bf LEMMA} \ %24
If \mT is a \bk-torus split by a finite Galois extension $\K/\bk$, then 
\[
X(T_\K) \ = \ \Mor_\K(T_\K,G_{m,\K})
\]
is a $\Gal(\K/\bk)$-lattice.
\end{x}
\vspace{0.1cm}

\begin{x}{\small\bf THEOREM} \ %25
Fix a finite Galois extension $\K/\bk$ $-$then the functor 
\[
T \ra X(T_\K)
\]
from the category of \bk-tori split by $\K/\bk$ to the category of $\Gal(\K/\bk)$-lattices is a contravariant equivalence of categories.
\end{x}
\vspace{0.1cm}

\begin{x}{\small\bf \un{N.B.}} \ %26
Suppose that \mT is a \bk-torus split by a finite Galois extension
%%----------------------------------------------------------------------------------------------06
$\K/\bk$.  
Form $\K[X(T_\K)]$, thus operationally, $\forall \ \sigma \in \Gal(\K/\bk)$, 
\[
\sigma \bigl( \sum\limits_i a_i \chi_i \bigr) \ = \ 
\sum\limits_i \sigma(a_i) \sigma(\chi_i) \qquad (a_i \in \K, \ \chi_i \in X(T_\K)).
\]
Pass now to the invariants 
\[
\K[X(T_\K)] \qquad (G = \Gal(\K/\bk)).
\]
Then 
\[
T \ \approx \  \Spec(\K[X(T_\K)]^G).
\]
And
\begin{align*}
T(A \ots{\bk} \K)^G \ 
&=\ T(A)\\
&\approx\ \Mor(\Spec(A),T)\\
&\approx\ \Mor(\Spec(A),\Spec(\K[X(T_\K)]^G)\\
&\approx\ \Mor_\bk(\K[X(T_\K)]^G,A)\\
&\approx\ \Mor_\K(\K[X(T_\K)],A \ots{\bk} \K)^G\\
&\approx\ \Mor_\Z(X(T_\K),(A \ots{\bk} \K)^\times)^G\\
&\approx\ \Mor_{\Z[G]} (X(T_\K),(A \ots{\bk} \K)^\times).
\end{align*}

[Note: \ 
Let $T = \Res_{\K/\bk} (G_{m,\K})$ $-$then on the one hand, 
\[
\Mor_{\Z[G]} (\Z[G], (A \ots{\bk} \K)^\times ) \ \approx \ (A \ots{\bk} \K)^\times,
\]
while on the other, 
\allowdisplaybreaks
\begin{align*}
\Res_{\K/\bk} (G_{m,\K})(A) \ 
&=\ (A \ots{\bk} \K)^\times
\\[12pt]
&\approx\ 
\Mor_{\Z[G]} (X(T_\K),(A \ots{\bk} \K)^\times).
\end{align*}
Therefore 
\[
X(T_\K) \ \approx  \ \Z[G].]
\]
\end{x}
\vspace{0.3cm}

%%----------------------------------------------------------------------------------------------07

Take $\bk = \R$, $\K = \C$, and let $\sigma$ be the nontrivial element of $\Gal(\C/\R)$ $-$then every 
$\R$-torus \mT gives rise to a free $\Z$-module of finite rank supplied with an involution corresponding to $\sigma$. 
And conversely \ldots \ .
\vspace{0.2cm}

There are three ``basic'' $\R$-tori.
\vspace{0.2cm}

\qquad\qquad 1. \quad $T = G_{m,\R}$.   \  In this case, 
\[
X(T_\C) \ = \ X(G_{m,\C}) \ \approx \  \Z
\]
and the Galois action is trivial.
\vspace{0.2cm}

\qquad\qquad 2. \quad $T = \Res_{\C/\R}(G_{m,\C})$.   \  In this case, 

\begin{align*}
X(T_\C) \ 
&\approx \ X(G_{m,\C} \times G_{m,\C}) \qquad (\text{cf.} \ \S3, \ \#26)\\
&\approx \ \Z \times \Z
\end{align*}
and the Galois action swaps coordinates.
\vspace{0.2cm}

\qquad\qquad 3. \quad $T = \SO_2$.  \  In this case, 
\begin{align*}
X((\SO_2)_\C) \ 
&\approx \ X(G_{m,\C}) \\
&\approx \  \Z
\end{align*}
and the Galois action is multiplication by $-1$.

[Note: \ 
\[
\SO_2:\bALG / \R \ra \bGRP
\]
is the functor defined by the rule
\[
\SO_2(A) \ = \ \bigl\{
\begin{pmatrix}
a &b\\
-b &a\\
\end{pmatrix}
: a, b \in A \ \& \  a^2 + b^2 = 1\bigr\}.
\]
Then $\SO_2$ is an algebraic $\R$-group such that 
\[
(\SO_2)_\C \ \approx \ G_{m,\C},
\]
%%----------------------------------------------------------------------------------------------08
so $\SO_2$ is an $\R$-torus and $\SO_2(\R)$ can be identified with 
$S$ ($= \{z \in \C: z\ov{z} = 1\}$).]
\vspace{0.3cm}

\begin{x}{\small\bf THEOREM} \ %27
Every $\R$-torus is isomorphic to a finite product of copies of the three basic tori described above.
\end{x}
\vspace{0.3cm}

Here is the procedure.  \ 
Fix a $\Z$-free module \mM of finite rank and an involution $\iota:M \ra M$ $-$then \mM can be decomposed as a direct sum 
\[
M_+ \oplus M_\sw \oplus M_-,
\]
where $\iota = 1$ on $M_+$, $\iota$ is a sum of 2-dimensional swaps on $M_\sw$ \ (or still, \ 
$M_\sw$ $=$ $\oplus \hsp \Z[\Gal(\C/\R)]$), and $\iota = -1$ on $M_-$.

\begin{x}{\small\bf SCHOLIUM} \ %28
If \mT is an $\R$-torus, then there exist unique nonnegative integers $a, b, c$ such that
\[
T(\R) \ \approx \ (\R^\times)^a \times (\C^\times)^b \times S^c.
\]
\end{x}
\vspace{0.3cm}

\begin{x}{\small\bf REMARK} \ %29
The classification of $\C$-tori is trivial:  Any such is a finite product of the $G_{m,\C}$.
\end{x}
\vspace{0.3cm}

\begin{x}{\small\bf RAPPEL} \ %30
Let $\K/\bk$ be a finite Galois extension and let \mA be a \bk-algebra $-$then there is a norm map
\[
 (A \ots{\bk} \K)^\times \ra A^\times \qquad (\approx (A \ots{\bk} \bk)^\times).
\]
\end{x}
\vspace{0.3cm}

\begin{x}{\small\bf CONSTRUCTION} \ %31
Let $\K/\bk$ be a finite Galois extension $-$then there is a norm map
\[
\tN_{\K/\bk}:\Res_{\K/\bk} (G_{m,\K}) \ra G_{m,\bk}.
\]
\vspace{0.2cm}
%%----------------------------------------------------------------------------------------------09

[For any \bk-algebra \mA, 
\begin{align*}
\Res_{\K/\bk}(G_{m,\K})(A) \ 
&=\ G_{m, \K}(A \ots{\bk} \K)\\
&=\ (A \ots{\bk} \K)^\times \ra A^\times \ = \ G_{m,\bk}(A).]
\end{align*}

[Note: \ $\tN_{\K/\bk}$ is not to be confused with the arrow of adjunction 
\[
G_{m,\bk} \ra \Res_{\K/\bk}(G_{m,\K}).]
\]
\end{x}
\vspace{0.3cm}

\begin{x}{\small\bf \un{N.B.}} \ %32
\[
\tN_{\K/\bk} \in X(\Res_{\K/\bk} (G_{m,\K})).
\]
\end{x}
\vspace{0.3cm}

\begin{x}{\small\bf NOTATION} \ %33
Let $\Res_{\K/\bk}^{(1)} (G_{m,\K})$ be the kernel of $\tN_{\K/\bk}$.
\end{x}
\vspace{0.3cm}

\begin{x}{\small\bf LEMMA} \ %34
$\Res_{\K/\bk}^{(1)} (G_{m,\K})$ is a \bk-torus and there is a short exact sequence
\[
1 \ra \Res_{\K/\bk}^{(1)} (G_{m,\K}) \ra \Res_{\K/\bk} (G_{m,\K}) \ra G_{m,\bk} \ra 1.
\]
\end{x}
\vspace{0.1cm}

\begin{x}{\small\bf EXAMPLE} \ %35
Take $\bk = \R$, $\K = \C$ $-$then 
\[
\Res_{\C/\R}^{(1)} (G_{m,\C}) \ \approx \ \SO_2
\]
and there is a short exact sequence
\[
1 \ra \SO_2 \ra \Res_{\C/\R} (G_{m,\C}) \ra G_{m,\R} \ra 1.
\]

[Note: \ On $\R$-points, this becomes
\[
1 \ra \tS \ra \C^\times \ra \R^\times \ra 1.]
\] 
\end{x}
\vspace{0.3cm}
%%----------------------------------------------------------------------------------------------10

\begin{x}{\small\bf DEFINITION} \ %36
Let \mT be a \bk-torus $-$then \mT is 
\un{\bk-anisotropic}
\index{\bk-anisotropic} 
if $X(T) = \{0\}$.
\end{x}
\vspace{0.1cm}

\begin{x}{\small\bf EXAMPLE} \ %37
$\SO_2$ is $\R$-anisotropic.
\end{x}
\vspace{0.3cm}

\begin{x}{\small\bf THEOREM} \ %38
Every \bk-torus \mT has a unique maxmal \bk-split subtorus $T_s$ and a unique maximal $\bk$-anisotropic subtorus $T_a$.  \ 
The intersection $T_s \hsp \cap \hsp \hsp T_a$ is finite and $T_s \hsp \cdot \hsp T_a = T$.
\end{x}
\vspace{0.3cm}

\begin{x}{\small\bf LEMMA} \ \  %39
$\Res_{\K/\bk}^{(1)} (G_{m,\K})$ is \bk-anisotropic.
\vspace{0.2cm}

PROOF \ 
Setting $G = \Gal(\K/\bk)$, under the functoriality of \#25, the norm map
\[
\tN_{\K/\bk}:\Res_{\K/\bk} (G_{m,\K}) \ra G_{m,\bk}
\]
corresponds to the homomorphism $\Z \ra \Z[G]$ of \mG-modules that sends $n$ to 
$n \hsp \bigl( \sum\limits_G \sigma \bigr)$, the quotient 
$\Z[G] / \Z\bigl( \sum\limits_G \sigma \bigr)$ being $X(T_\K)$, where
\[
T \ = \ \Res_{\K/\bk}^{(1)} (G_{m,\K}).
\]
And
\[
\Z[G]^G 
\ = \ 
\Z\bigl( \sum\limits_G \sigma \bigr).
\]
\end{x}
\vspace{0.3cm}

\begin{x}{\small\bf \un{N.B.}} \ \ %40 
$\Res_{\K/\bk}^{(1)} (G_{m,\K})$ is the maximal \bk-anisotropic subtorus of $\Res_{\K/\bk} (G_{m,\K})$.
\end{x}
\vspace{0.3cm}

\begin{x}{\small\bf DEFINITION} \ %41
Let \mG, \mH be algebraic \bk-groups $-$then a homomorphism $\phi:G \ra H$ is an \un{isogeny} if it is surjective with a finite kernel.
\end{x}
\vspace{0.3cm}

%%----------------------------------------------------------------------------------------------11
\begin{x}{\small\bf DEFINITION} \ %42
Let \mG, \mH be algebraic \bk-groups $-$then \mG, \mH are said to be 
\un{isogeneous}
\index{isogeneous} 
if there is an isogeny between them.
\end{x}
\vspace{0.3cm}

\begin{x}{\small\bf THEOREM} \ %43
Two \bk-tori $T^\prime$, $T\pp$ per \#25 are isogeneous iff the $\Q[\Gal(\K / \bk)]$-modules
\[
\begin{cases}
\ X(T_\K^\prime) \ots{\Z} \Q \\
\ X(T_\K\pp) \ots{\Z} \Q 
\end{cases}
\]
are isomorphic.
\end{x}
\vspace{0.1cm}
%%%%%%%%%%%%%%%%%%%%%%%%%%%%%%%%%%%%%%
%%%%%%%%%%%%%%%%%%%%%%%%%%%%%%%%%%%%%%
%%%%%%%%%%%%%%%%%%%%%%%%%%%%%%%%%%%%%%

%% file: _05_THE_LLC.tex
\chapter{
$\boldsymbol{\S}$\textbf{5}.\quad  THE LLC}
\setlength\parindent{2em}
\setcounter{theoremn}{0}
%%----------------------------------------------------------------------------------------------01

\ \indent 

\begin{x}{\small\bf \un{N.B.}}  \ %01
The term ``LLC'' means ``local Langlands correspondence'' (cf. \#26).  

Let $\K$ be a non-archimedean local field $-$then the image of $\rec_\K:\K^\times \ra G_\K^{\ab}$ 
is $W_\K^\ab$ and the induced map $\K^\times \ra W_\K^\ab$ is an isomorphism of topological groups.
\end{x}
\vspace{0.1cm}

\begin{x}{\small\bf SCHOLIUM} \ %02
There is a bijective correspondece between the characters of $W_\K$ and the characters of $\K^\times$: 
\[
\Mor(W_\K,\C^\times) \ \approx \ \Mor(\K^\times,\C^\times).
\]

[Note: \ 
``Character'' means continuous homomorphism.  
So, if $\chi:W_\K \ra \C^\times$ is a character, then $\chi$ must be trivial on $W_\K^*$ 
($\C^\times$ being abelian), hence by continuity, trivial on $\ov{W_\K^*}$, thus $\chi$ factors through 
$W_\K / \hsp \ov{W_\K^*} = W_\K^\ab$.]

Let \mT be a $\K$-torus $-$then \mT is isomorphic to a closed subgroup of some 
$\GL_{n,\K}$ $(\exists \ n)$.  
But $\GL_{n,\K}(\K)$ is a locally compact topological group, thus $T(\K)$ is a locally compact topological group 
(which, moreover, is abelian).
\end{x}
\vspace{0.1cm}

\begin{x}{\small\bf \un{N.B.}}  \ %03
For the record, 
\[
G_{m,\K}(\K) \ = \ \K^\times \ = \ \GL_{1,\K}(\K).
\]
\end{x}
\vspace{0.1cm}

\begin{x}{\small\bf EXAMPLE} \ %04
Let $\LL/\K$ be a finite extension and consider $T = \Res_{\LL/\K}(G_{m,\LL})$ $-$then $T(\K) = \LL^\times$.
\end{x}
\vspace{0.1cm}

Roughly speaking, the objective now is to describe $\Mor(T(\K),\C^\times)$ in terms of data attached to $W_\K$ 
but to even state the result requires some preparation.

\vspace{0.2cm}

%%----------------------------------------------------------------------------------------------02

\begin{x}{\small\bf \un{N.B.}}  \ %05
The case when $T = G_{m,\K}$ is local class field theory \ldots \ .
\end{x}
\vspace{0.1cm}

\begin{x}{\small\bf EXAMPLE} \ %06
Suppose that \mT is $\K$-split:
\[
T \ \approx \ G_{m,\K}\times \cdots \times G_{m,\K} \qquad (\text{$d$ factors}).
\]
Then
\begin{align*}
\prod\limits_{i = 1}^d \Mor(W_\K,\C^\times) \ 
&\approx\ \prod\limits_{i = 1}^d \Mor(\K^\times,\C^\times)\\
&\approx\ \Mor\bigl(\prod\limits_{i = 1}^d  \K^\times,\C^\times \bigr)\\
&\approx\ \Mor(T(\K),\C^\times).
\end{align*}

Given a $\K$-torus \mT, put
\[
\begin{cases}
\ X^*(T) \ = \ \Mor_{\K^\sep}(T_{\K^\sep}, G_{m,\K^\sep})\\
\ X_*(T) \ = \ \Mor_{\K^\sep}(G_{m,\K^\sep},T_{\K^\sep})
\end{cases}
.
\]
\end{x}
\vspace{0.1cm}

\begin{x}{\small\bf LEMMA} \ %07
Canonically, 
\[
 X_*(T) \ots{\Z} \C^\times \ \approx \ \Mor(X^*(T),\C^\times).
\]

\vspace{0.1cm}

PROOF \ 
Bearing in mind that 
\[
\Mor_{\K^\sep}(G_{m,\K^\sep},G_{m,\K^\sep})  \ \approx \ \Z,
\]
define a pairing 
\[
\begin{tikzcd}[sep=large]
{X^*(T) \times X_*(T)} 
\ar{r}{\langle \ , \ \rangle} 
&{\Z}
\end{tikzcd}
\]
%%----------------------------------------------------------------------------------------------03
by sending $(\chi^*,\chi_*)$ to $\chi^* \circ \chi_* \in \Z$.  
This done, given $\chi_* \otx z$, assign to it the homomorphism 
\[
\chi^* \ \ra z^{\langle \chi^*,\chi_* \rangle}.
\]
\end{x}
\vspace{0.1cm}

\begin{x}{\small\bf NOTATION} \ %08
Given a $\K$-torus \mT, put
\[
\widehat{T} \ = \ \Spec(\C[X_*(T)]).
\]
\end{x}
\vspace{0.1cm}

\begin{x}{\small\bf LEMMA} \ %09
$\widehat{T}$ is a split $\C$-torus such that 
\[
\begin{cases}
\ X^*(\widehat{T}) \ \equiv \  \Mor_\C(\widehat{T},G_{m,\C}) \ \approx \ X_*(T)\\
\ X_*(\widehat{T}) \ \equiv \  \Mor_\C(G_{m,\C},\widehat{T}) \ \approx \ X^*(T)
\end{cases}
.
\]
Therefore
\begin{align*}
\Mor(X_*(T),\C^\times) \ 
&\approx\ \Mor(X^*(\widehat{T}),\C^\times)\\
&\approx\ X_*(\widehat{T}) \ots{\Z} \C^\times\\
&\approx\ X^*(T) \ots{\Z} \C^\times.
\end{align*}
\end{x}
\vspace{0.1cm}

\begin{x}{\small\bf LEMMA} \ %10
\[
\widehat{T}(\C) \ \approx \ X^*(T) \ots{\Z} \C^\times.
\]

\vspace{0.1cm}

PROOF \ 
In fact, 
\begin{align*}
\widehat{T}(\C) \ 
&\approx\ \Mor(X^*(\widehat{T}),\C^\times) \qquad (\text{cf.} \ \S4, \ \#12)\\
&\approx\ \Mor(X_*(T),\C^\times)\\
&\approx\ X^*(T) \ots{\Z} \C^\times.
\end{align*}
\end{x}
\vspace{0.1cm}

%%----------------------------------------------------------------------------------------------04

\begin{x}{\small\bf DEFINITION} \ %11
$\widehat{T}$ is the 
\un{complex dual torus}
\index{complex dual torus} 
of \mT.
\end{x}
\vspace{0.1cm}

\begin{x}{\small\bf EXAMPLE} \ %12
Under the assumptions of \#6, 
\begin{align*}
\widehat{T}(\C) \ 
&\approx\ X^*(T) \ots{\Z} \C^\times\\
&\approx\ \Z^d \ots{\Z} \C\\
&\approx(\C^\times)^d.
\end{align*}
Therefore
\begin{align*}
\Mor(W_\K,\widehat{T}(\C)) \ 
&\approx\ \Mor(W_\K,(\C^\times)^d)\\
&\approx\ \prod\limits_{i = 1}^d \Mor(W_\K,\C^\times)\\
&\approx\ \Mor(T(\K),\C^\times).
\end{align*}
\end{x}
\vspace{0.1cm}

\begin{x}{\small\bf RAPPEL} \ %13
If \mG is a group and if \mA is a \mG-module, then 
\[
H^1(G,A) \ = \ \frac{Z^1(G,A)}{B^1(G,A)}.
\]
\end{x}
\vspace{0.1cm}

\qquad\qquad \textbullet \quad 
$Z^1(G,A)$ (the 
\un{1-cocycles}) 
\index{1-cocycles}
consists of those maps $f:G \ra A$ such that $\forall \ \sigma, \tau \in G$, 
\[
f(\sigma \tau) \ = \ f(\sigma) + \sigma (f(\tau)).
\]

\qquad\qquad \textbullet \quad
$B^1(G,A)$  (the 
\un{1-coboundaries})
\index{1-coboundaries} 
consists of those maps $f:G \ra A$ for which $\exists$ an $a \in A$ such that 
$\forall \ \sigma \in G$, 
\[
f(\sigma) \ = \ \sigma a - a.
\]

[Note: \ 
\[
H^1(G,A) \ = \ \Mor(G,A)
\]
if the action is trivial.]
\vspace{0.3cm}
%%----------------------------------------------------------------------------------------------05

\begin{x}{\small\bf NOTATION} \ %14
If \mG is a topological group and if \mA is a topological \mG-module, then
\[
\Mor_c(G,A)
\]
is the group of continuous group homomorphisms from \mG to \mA.  
Analogously, 
\[
\begin{cases}
\ Z_c^1(G,A) = \text{``continuous 1-cocycles''}\\
\ B_c^1(G,A) = \text{``continuous 1-coboundaries''}
\end{cases}
\]
and 
\[
H_c^1(G,A) \ = \ \frac{Z_c^1(G,A)}{B_c^1(G,A)}.
\]

Let \mT be a $\K$-torus $-$then $G_\K$ ($= \Gal(\K^\sep/\K)$) operates on $X^*(G)$, thus $W_\K \subset G_\K$ 
operates on $X^*(G)$ by restriction.  
Therefore $\widehat{T}(\C)$ is a $W_\K$-module, so it makes sense to form 
\[
H_c^1(W_\K,\widehat{T}(\C)).
\]
\end{x}
\vspace{0.1cm}

\begin{x}{\small\bf NOTATION} \ %15
$\bTOR_\K$ is the category of $\K$-tori.
\end{x}
\vspace{0.1cm}

\begin{x}{\small\bf LEMMA} \ %16
The assignment 
\[
T \ra H_c^1(W_\K,\widehat{T}(\C))
\]
defines a functor 
\[
\bTOR_\K^\OP \lra \bAB.
\]

%%----------------------------------------------------------------------------------------------06

[Note: \ 
Suppose that $T_1 \ra T_2$ $-$then 
\[
(T_1)_{\K^\sep} \ra (T_2)_{\K^\sep}
\]
\qquad\qquad $\implies$
\[
X^*(T_2) \ra X^*(T_1)
\]
\qquad\qquad $\implies$
\[
\widehat{T}_2(\C) \ra \widehat{T}_1(\C)
\]
\qquad\qquad $\implies$
\[
H_c^1(W_\K,\widehat{T}_2(\C)) \ \ra \ H_c^1(W_\K,\widehat{T}_1(\C)).]
\]
\end{x}
\vspace{0.1cm}

\begin{x}{\small\bf LEMMA} \ %17
The assignment 
\[
T \ra \Mor_c(T(\K),\C^\times)
\]
defines a functor 
\[
\bTOR_\K^\OP \lra \bAB.
\]
\end{x}
\vspace{0.1cm}

\begin{x}{\small\bf THEOREM} \ %18
The functors 
\[
T \ra H_c^1(W_\K,\widehat{T}(\C))
\]
and 
\[
T \ra \Mor_c(T(\K),\C^\times)
\]
are naturally isomorphic.
\end{x}
\vspace{0.1cm}

\begin{x}{\small\bf SCHOLIUM} \ %19
There exist isomorphisms
\[
\iota_T: H_c^1(W_\K,\widehat{T}(\C)) \ra \Mor_c(T(\K),\C^\times)
\]
%%----------------------------------------------------------------------------------------------07
such that if $T_1 \ra T_2$, then the diagram 
\[
\begin{tikzcd}[sep=large]
{H_c^1(W_\K,\widehat{T}_1(\C))} \ar{rr}{\iota_{T_1}}
&&{\Mor_c(T_1(\K),\C^\times)}\\
{H_c^1(W_\K,\widehat{T}_2(\C))} \ar{u}  \ar{rr}[swap]{\iota_{T_2}}
&&{\Mor_c(T_2(\K),\C^\times)} \ar{u}
\end{tikzcd}
\]
commutes.

\end{x}
\vspace{0.1cm}

\begin{x}{\small\bf EXAMPLE} \ %20
Under the assumptions of \#12, the action of $G_\K$ is trivial, hence the action of $W_\K$ is trivial.  Therefore
\begin{align*}
H_c^1(W_\K,\widehat{T}(\C)) \
&=\ \Mor_c(W_\K,\widehat{T}(\C))
\\[12pt]
&\approx\ \Mor_c(T(\K),\C^\times).
\end{align*}

[Note: \ 
The earlier use of the symbol Mor tacitly incorporated ``continuity''.]
\end{x}
\vspace{0.1cm}

There is a special case that can be dealt with directly, viz. when $\LL/\K$ is a finite Galois extension and 
\[
T \ = \ \Res_{\LL/\K} (G_{m,\LL}).
\]
The discussion requires some elementary cohomological generalities which have been collected in the Appendix below.

\vspace{0.2cm}

\begin{x}{\small\bf RAPPEL} \ %21
$W_\LL$ is a normal subgroup of $W_\K$ of finite index:
\[
W_\K/W_\LL \ \approx G_\K / G_\LL \ \approx \ \Gal(\LL/\K).
\]
Proceeding, 
\[
T_{\K^\sep} \ \approx \ \prod\limits_{\sigma \in \Gal(\LL/\K)} \sigma G_{m,\LL} \qquad (\text{cf.} \ \#6),
\]
%%----------------------------------------------------------------------------------------------08
so 
\[
X^*(T) \ \approx \ \Z[W_\K/W_\LL],
\]
where
\begin{align*}
\Z[W_\K/W_\LL] \
&\approx\ \Ind_{W_\LL}^{W_{\K_\Z}} \\
&\equiv\ \Z[W_\K] \ots{\Z[W_\LL]} \Z.
\end{align*}
It therefore follows that 
\begin{align*}
\widehat{T}(\C) \ 
&\approx\ X^*(T) \ots{\Z} \C^\times
\\[12pt]
&\approx\ \Z[W_\K] \ots{\Z[W_\LL]}\Z \ots{\Z} \C^\times
\\[12pt]
&\approx\ \Z[W_\K] \ots{\Z[W_\LL]} \C^\times 
\\[12pt]
&\equiv\ \Ind_{W_\LL}^{W_{\K_{\C^\times}}}.
\end{align*}
Consequently
\begin{align*}
H^1(W_\K,\widehat{T}(\C)) \ 
&\approx\ H^1(W_\K,\Ind_{W_\LL}^{W_{\K_{\C^\times}}})
\\[12pt]
&\approx\ H^1(W_\LL,\C^\times) \qquad\qquad (\text{Shapiro's lemma})
\\[12pt]
&\approx\ \Mor(W_\LL,\C^\times)
\\[12pt]
&\approx\ \Mor(\LL^\times,\C^\times)
\\[12pt]
&\approx\ \Mor(T(\K),\C^\times),
\end{align*}
which completes the proof modulo ``continuity details'' that we shall not stop to sort out.
\end{x}
\vspace{0.1cm}

%%----------------------------------------------------------------------------------------------09

\begin{x}{\small\bf DEFINITION} \ %22
The 
\un{L-group}
\index{L-group} 
of \mT is the semidirect product
\[
\text{\mLT} \ = \ \widehat{T}(\C) \rtimes W_\K.
\]

Because of this, it will be best to first recall ``semidirect product theory''.
\end{x}
\vspace{0.1cm}

\begin{x}{\small\bf RAPPEL} \ %23
If \mG is a group and if \mA is a \mG-module, then there is a canonical extension of \mG by \mA, namely
\[
0 \ra A \overset{i}{\lra} A \rtimes G \overset{\pi}{\lra} G \ra 1,
\]
where $A \rtimes G$ is the semidirect product.
\end{x}
\vspace{0.1cm}

\begin{x}{\small\bf DEFINITION} \ %24
A 
\un{splitting}
\index{splitting} 
of the extension 
\[
0 \ra A \overset{i}{\lra} A \rtimes G \overset{\pi}{\lra} G \ra 1
\]
is a homomorphism $s:G \ra A \rtimes G$ such that $\pi \circ s = \id_G$.
\end{x}
\vspace{0.1cm}

\begin{x}{\small\bf FACT} \ %25
The splittings of the extension
\[
0 \ra A \overset{i}{\lra} A \rtimes G \overset{\pi}{\lra} G \ra 1
\]
determine and are determined by the elements of $\Z^1(G,A)$.
\end{x}
\vspace{0.1cm}

Two splittings $s_1$, $s_2$ are said to be 
\un{equivalent}
\index{splittings \\ equivalent} 
if there is an element $a \in A$ such that 
\[
s_1(\sigma) \ = \ i(a) s_2(\sigma) i(a)^{-1} \quad (\sigma \in G).
\]
If
\[
\begin{cases}
\ f_1 \longleftrightarrow s_1\\
\ f_2 \longleftrightarrow s_2
\end{cases}
\]
%%----------------------------------------------------------------------------------------------10
are the 1-cocycles corresponding to 
$
\begin{cases}
\ s_1\\
\ s_2\\
\end{cases}
$
, then their difference $f_2 - f_1$ is a 1-coboundary.
\vspace{0.2cm}

\begin{x}{\small\bf SCHOLIUM} \ %26
The equivalence classes of splittings of the extension 
\[
0 \ra A \overset{i}{\lra} A \rtimes G \overset{\pi}{\lra} G \ra 1
\]
are in a bijective correspondence with the elements of $H^1(G,A)$.
\end{x}
\vspace{0.1cm}

Return now to the extension 
\[
\begin{tikzcd}%[sep=large]
{0 \lra \widehat{T}(\C) \lra \widehat{T}(\C) \rtimes W_\K \lra W_\K \lra 1}
\quad \arrow[d,shift right=0.5,dash] \arrow[d,shift right=-0.5,dash] \\
{\text{\mLT}}
\end{tikzcd}
\]
but to reflect the underlying topologies, work with continuous splittings and call them 
\un{admissible homomorphisms}.  
\index{admissible homomorphisms}
Introducing the obvious notion of equivalence, denote by $\Phi_\K(T)$ the set of equivalence classes of admissible homomorphisms, hence 
\[
\Phi_\K(T) \ \approx \ H_c^1(W_\K,\widehat{T}(\C)).
\]
On the other hand, denote by $\sA_\K(T)$ the group of characters of $T(\K)$, i.e., 
\[
\sA_{\K}(T)  \ \approx \ \Mor_c(T(\K),\C^\times).
\]

\vspace{0.2cm}

\begin{x}{\small\bf THEOREM} \ %27
There is a canonical isomorphism
\[
\Phi_\K(T) \ra \sA_{\K}(T).
\]

[This statement is just a rephrasing of \#18 and is the LLC for tori.]
\end{x}
\vspace{0.1cm}

\begin{x}{\small\bf HEURISTICES} \ %28
To each admissible homomorphism of $W_\K$ into \mLT, it is possible to associate an irreducible automorphic representation of 
$T(\K)$ (a.k.a. a character of $T(\K)$)  and all such arise in this fashion.
\end{x}
\vspace{0.1cm}
%%----------------------------------------------------------------------------------------------11

It remains to consider the archimedean case: $\C$ or $\R$.

\vspace{0.1cm}

\qquad\qquad \textbullet \quad If \mT is a $\C$-torus, then \mT is isomorphic to a finite product
\[
G_{m,\C} \times \cdots \times G_{m,\C}
\]
and 
\begin{align*}
T(\C) \ 
&\approx\ \Mor(X^*(T),\C^\times)
\\[12pt]
&\approx\ X_*(T) \ots{\Z} \C^\times.
\end{align*}
Furthermore, $W_\C = \C^\times$ and the claim is that
\[
H_c^1(W_\C,\widehat{T}(\C)) \ \equiv \ \Mor_c(\C^\times,\widehat{T}(\C))
\]
is isomorphic to 
\[
\Mor_c(T(\C),\C^\times).
\]
But
\begin{align*}
\Mor_c(\C^\times,\widehat{T}(\C)) \ 
&\approx\ \Mor_c(\C^\times,X^*(T) \ots{\Z} \C^\times)\\
&\approx\ \Mor_c(\C^\times, \Mor(X_*(T),\C^\times))\\
&\approx\ \Mor_c(X_*(T) \ots{\Z} \C^\times,\C^\times)\\
&\approx\ \Mor_c(T(\C),\C^\times).
\end{align*}

\qquad\qquad \textbullet \quad  If \mT is a $\R$-torus, then \mT is isomorphic to a finite product
\[
(G_{m,\R})^a \times (\Res_{\C/\R}(G_{m,\C}))^b \times (\SO_2)^c
\]
and it is enough to look at the three irreducible possibilities.

%%----------------------------------------------------------------------------------------------12

1. \quad $T = G_{m,\R}$. \ The point here is that $W_\R^{ab} \approx \R^\times \equiv T(\R)$.

2. \quad $T = \Res_{\C/\R}(G_{m,\C})$. \ One can imitate the argument used above for its non-archimedean analog.

3. \quad $T = \SO_2$. 
The initial observation is that $X(T) = \Z$ with action $n \ra -n$, so $\widehat{T}(\C) = \C^\times$ with action $z \ra \ds\frac{1}{z}$.  
And \ldots .
\vspace{0.5cm}

\[
\textbf{APPENDIX}
\]
\setcounter{theoremn}{0}
\vspace{0.5cm}

Let \mG be a group (written multiplicatively).
\vspace{0.3cm}

\begin{x}{\small\bf DEFINITION} \ %01
A left (right) \mG-module is an abelian group \mA equipped with a left (right) action of \mG, 
i.e., with a homomorphism $G \ra \Aut(A)$.
\end{x}
\vspace{0.1cm}

\begin{x}{\small\bf \un{N.B.}}  \ %02
Spelled out, to say that \mA is a left \mG-module means that there is a map
\[
\begin{cases}
\ \mG \times A \ra A\\
\ (\sigma,a) \ra \sigma a
\end{cases}
\]
such that 
\[
\tau(\sigma a) \ = \ (\tau \sigma) a, \qquad 1 a \ = \ a,
\]
thus \mA is first of all a left \mG-set.   
To say that \mA is a left \mG-module then means in addition that
\[
\sigma (a + b) \ = \ \sigma a + \sigma b.
\]

[Note: \ 
For the most part, the formalities are worked out from the left, the agreement being that 
\[
\text{``left \mG-module'' \ = \ ``\mG-module''.]}
\]
\end{x}
\vspace{0.1cm}

%%----------------------------------------------------------------------------------------------13

\begin{x}{\small\bf NOTATION} \ %03
The group ring $\Z[G]$ is the ring whose additive group is the free abelian group with basis \mG and whose multiplication 
is determined by the multiplication in \mG and the distributive law.
\end{x}
\vspace{0.1cm}

A typical element of $\Z[G]$ is 
\[
\sum\limits_{\sigma \in G} m_\sigma \sigma,
\]
where $m_\sigma \in \Z$ and $m_\sigma = 0$ for all but finitely many $\sigma$.

\vspace{0.2cm}

\begin{x}{\small\bf \un{N.B.}}  \ %04
A \mG-module is the same thing as a $\Z[G]$-module.
\end{x}

\begin{x}{\small\bf LEMMA} \ %05
Given a ring \mR, there is a canonical bijection 
\[
\Mor(\Z[G],R) \ \approx \ \Mor(G,R^\times).
\]
\end{x}
\vspace{0.1cm}

\begin{x}{\small\bf CONSTRUCTION} \ %06
Given a \mG-set \mX, form the free abelian group $\Z[X]$ generated by \mX and extend the action of \mG on \mX 
to a $\Z$-linear action of \mG on $\Z[X]$ $-$then the resulting \mG-module is called a 
\un{permutation module}.
\index{permutation module}
\end{x}
\vspace{0.1cm}

\begin{x}{\small\bf EXAMPLE} \ %07
Let \mH be a subgroup of \mG and take $X = G/H$ (here \mG operates on $G/H$ by left translation), from which 
$\Z[G/H]$.
\end{x}
\vspace{0.1cm}

\begin{x}{\small\bf DEFINITION} \ %08
A 
\un{\mG-module homomorphism}
\index{\mG-module homomorphism} 
is a $\Z[G]$-module homomorphism.
\end{x}
\vspace{0.1cm}

\begin{x}{\small\bf NOTATION} \ %09
$\bMOD_G$ is the category of \mG-modules.
\end{x}
\vspace{0.1cm}

\begin{x}{\small\bf NOTATION} \ %10
Given \mA, \mB in $\bMOD_G$, write $\Hom_G(A,B)$ in place of $\Mor(A,B)$.
\end{x}
\vspace{0.1cm}

\begin{x}{\small\bf LEMMA} \ %11
Let $A, B  \in \bMOD_G$ $-$then $A \ots{\Z} B$ carries the \mG-module structure 
%%----------------------------------------------------------------------------------------------14
defined by $\sigma (a \otx a^\prime) = \sigma a \otx \sigma a^\prime$ and $\Hom_\Z(A,B)$ carries the \mG-module structure 
defined by $(\sigma \phi)(a) = \sigma\phi(\sigma^{-1} a)$.
\end{x}
\vspace{0.1cm}

\begin{x}{\small\bf LEMMA} \ %
If $G^\prime$ is a subgroup of \mG, then there is a homomorphism $\Z[G^\prime] \ra \Z[G]$ of rings and a functor 
\[
\Res_{G^\prime}^G:\bMOD_G \ra \bMOD_{G^\prime}
\]
of restriction.
\end{x}
\vspace{0.1cm}

\begin{x}{\small\bf DEFINITION} \ %13
Let $G^\prime$ be a subgroup of \mG, $-$then the 
\un{functor of induction}
\index{functor of induction}
\[
\Ind_{G^\prime}^G:\bMOD_{G^\prime} \ra \bMOD_{G}
\]
sends $A^\prime$ to 
\[
\Z[G] \ots{\Z[G^\prime]} A^\prime.
\]

[Note: \ 
$Z[G]$ is a right $\Z[G^\prime]$-module and $A^\prime$ is a left $\Z[G^\prime]$-module.  
Therefore the tensor product
\[
\Z[G] \ots{\Z[G^\prime]} A^\prime 
\]
is an abelian group.  And it becomes a left \mG-module under the operation 
$\sigma(r \otx a^\prime) = \sigma r \otx a^\prime$.]
\end{x}
\vspace{0.1cm}

\begin{x}{\small\bf EXAMPLE} \ %14
Let \mH be a subgroup of \mG.  
Suppose that \mH operates trivially on $\Z$ $-$then 
\[
\Z[G/H] \ \approx \ \Ind_H^G \Z.
\]
\end{x}
\vspace{0.1cm}

\begin{x}{\small\bf FROBENIUS RECIPROCITY} \ %15
$\forall \ A$ in $\bMOD_G$, $\forall \ A^\prime$ in $\bMOD_{G^\prime}$,
\[
\Hom_{G^\prime}(A^\prime,\Res_{G^\prime}^G A) \ \approx \ 
\Hom_G(\Ind_{G^\prime}^G A^\prime, A).
\]
\end{x}
\vspace{0.1cm}

%%----------------------------------------------------------------------------------------------15

\begin{x}{\small\bf REMARK} \ %16
$\forall \ A$ in $\bMOD_G$,
\[
\Ind_{G^\prime}^G \circ \Res_{G^\prime}^G A \ \approx \ \Z[G/G^\prime] \ots{\Z[G]} A.
\]

[\mG operates on the right hand side diagonally: $\sigma(r \otx a) = \sigma r \otx \sigma a$.]
\end{x}
\vspace{0.1cm}

\begin{x}{\small\bf LEMMA} \ %17
There is an arrow of inclusion
\[
\Z[G] \ots{\Z[G^\prime]} A^\prime \ra \Hom_{G^\prime} (\Z[G],A^\prime)
\]
which is an isomorphism if $[G:G^\prime] < \infty$.
\end{x}
\vspace{0.1cm}

\begin{x}{\small\bf NOTATION} \ %18
Given a \mG-module \mA, put
\[
A^G \ =  \{a \in A: \sigma a = a \ \forall \ \sigma \in G\}.
\]

[Note: \ 
$A^G$ is a subgroup of \mA, termed the 
\un{invariants}
\index{invariants} 
in \mA.]
\end{x}
\vspace{0.1cm}

\begin{x}{\small\bf LEMMA} \ %19
$A^G = \Hom_G(\Z,A)$ (trivial \mG-action on $\Z$).

[Note: \ 
By comparison,
\[
A \ = \ \Hom_G(\Z[G],A).]
\]
\end{x}
\vspace{0.1cm}

\begin{x}{\small\bf LEMMA} \ %20
$
\ \Hom_\Z(A,B)^G \ = \ \Hom_G(A,B).
$
\end{x}
\vspace{0.1cm}

$\bMOD_G$ is an abelian category.  As such, it has enough injectives (i.e., every \mG-module can be embedded in an 
injective \mG-module).

\vspace{0.2cm}

\begin{x}{\small\bf DEFINITION} \ %21
The 
\un{group cohomology}
\index{group cohomology} 
functor 
$H^q(G,-):\bMOD_G \ra \bAB$ 
is the right derived functor of $(-)^G$.

\vspace{0.1cm}

[Note: \ Recall the procedure:  \ To compute $H^q(G,A)$, choose an injective
%%----------------------------------------------------------------------------------------------16
resolution
\[
0 \ra A \ra I^0 \ra I^1 \ra \cdots  \ .
\]
Then $H^*(G,A)$ is the cohomology of the complex $(I)^G$.  
In particular: $H^0(G,A) = A^G$.]
\end{x}
\vspace{0.1cm}

\begin{x}{\small\bf LEMMA} \ %22
$H^q(G,A)$ is independent of the choice of injective resolutions.
\end{x}
\vspace{0.1cm}

\begin{x}{\small\bf LEMMA} \ %23
$H^q(G,A)$ is a covariant functor of \mA.
\end{x}
\vspace{0.1cm}

\begin{x}{\small\bf LEMMA} \ %24
If
\[
0 \ra A \ra B \ra C \ra 0
\]
is a short exact sequence of \mG-modules, then there is a functorial long exact sequence
\begin{align*}
0 
&\ra H^0(G,A) \ra H^0(G,B) \ra H^0(G,C)\\
&\ra H^1(G,A) \ra H^1(G,B) \ra H^1(G,C) \ra H^2(G,A) \ra \cdots\\
\cdots 
&\ra H^q(G,A) \ra H^q(G,B) \ra H^q(G,C) \ra H^{q+1}(G,A) \ra \cdots\\
\end{align*}
in cohomology.
\end{x}
\vspace{0.1cm}

\begin{x}{\small\bf \un{N.B.}}  \ %25
If $G = \{1\}$ is the trivial group, then 
\[
H^0(G,A) \ = \ A, \quad H^q(G,A) = 0 \qquad (q > 0).
\]

\vspace{0.1cm}

[Note: \ 
Another point is that for any \mG, every injective \mG-module \mA is cohomologically acyclic:
\[
\forall \ q > 0, \  H^q(G,A) = 0.]
\]
\end{x}
\vspace{0.1cm}

\begin{x}{\small\bf THEOREM (SHAPIRO'S LEMMA)} \ %26
If $[G:G^\prime] < \infty$, then $\forall \ q$, 
\[
H^q(G^\prime,A^\prime) \ \approx \ H^q(G,\Ind_{G^\prime}^G A^\prime).
\]
\end{x}
\vspace{0.1cm}

%%----------------------------------------------------------------------------------------------17
\begin{x}{\small\bf EXAMPLE} \ %27
Take $A^\prime = \Z[G^\prime]$ $-$then 
\begin{align*}
H^q(G^\prime,\Z[G^\prime]) \ 
&\approx\ H^q(G,\Z[G] \ots{\Z[G^\prime]} \Z[G^\prime])
\\[12pt]
&\approx\ H^q(G,\Z[G]).
\end{align*}
\end{x}
\vspace{0.1cm}

\begin{x}{\small\bf EXAMPLE} \ %28
Take $G^\prime = \{1\}$ (so \mG is finite) $-$then $\Z[G^\prime] = \Z$ and 
\[
H^q(\{1\},\Z) \ \approx \ H^q(G,\Z[G]).
\]
But the LHS vanishes if $q > 0$, thus the same is true of the RHS.  However, this fails if \mG is infinite.  
E.g.: \ Take for \mG the infinite cyclic group: \ $H^1(G,\Z[G]) \approx \Z$.
\vspace{0.1cm}

[Note: \ 
If \mG  is finite,  then \ $H^0(G,\Z[G]) \approx  \Z$ \  while  if  \mG  is  infinite,  then $H^0(G,\Z[G]) = 0$.]
\end{x}
\vspace{0.1cm}

\begin{x}{\small\bf EXAMPLE} \ %29
Take $A^\prime = \Z$ $-$then 
\begin{align*}
H^q(G^\prime,\Z) \ 
&\approx\ H^q(G,\Ind_{G^\prime}^G \Z)
\\[12pt]
&\approx\ H^q(G,\Z[G/G^\prime]).
\end{align*}
\end{x}
\vspace{0.1cm}

%%%%%%%%%%%%%%%%%%%%%%%%%%%%%%%%%%%%%%
%%%%%%%%%%%%%%%%%%%%%%%%%%%%%%%%%%%%%%
%%%%%%%%%%%%%%%%%%%%%%%%%%%%%%%%%%%%%%

%% file: _06_TAMAGAWA_MEASURES.tex
\chapter{
$\boldsymbol{\S}$\textbf{6}.\quad  TAMAGAWA MEASURES}
\setlength\parindent{2em}
\setcounter{theoremn}{0}
%%----------------------------------------------------------------------------------------------01

\ \indent 

Suppose given a $\Q$-torus \mT of dimension $d$ $-$then one can introduce

\[
\begin{tikzcd}[sep=small]
{T(\Q) \subset T(\R), \quad T(\Q) \subset T(\Q_p)} \\
{\hspace{4.3cm}\bigcup}\\
{\hspace{4.3cm}T(\Z_p)}
\end{tikzcd}
\]
and 
\[
T(\Q) \subset T(\A).
\]

\vspace{0.2cm}

\begin{x}{\small\bf EXAMPLE} \ %01
Take $T = G_{m,\Q}$ $-$then the above data becomes 
\[
\begin{tikzcd}[sep=small]
{\Q^\times \subset \R^\times, \quad \Q^\times \subset \Q_p^\times} \\
{\hspace{3.1cm}\bigcup}\\
{\hspace{3.1cm}\Z_p^\times}
\end{tikzcd}
\]
and
\[
\Q^\times \subset \A^\times \ = \I.
\]

\end{x}
\vspace{0.3cm}

\begin{x}{\small\bf LEMMA} \ %02
$T(\Q)$ is a discrete subgroup of $T(\A)$.
\end{x}
\vspace{0.3cm}

\begin{x}{\small\bf RAPPEL} \ %03
$\I^1 = \Ker \acdot_\A$, where for $x \in \I$,
\[
\abs{x}_\A \ = \ \prod\limits_{p \leq \infty} \abs{x_p}_p.
\]
And the quotient $\I^1 / \Q^\times$ is a compact Hausdorff space.

Each $\chi \in X(T)$ generates continuous homomorphisms
\[
\begin{cases}
\begin{tikzcd}[sep=small]
{\chi_p:T(\Q_p) \ra \Q_p^\times} \ar{rr}{\acdot_p} &&{\R_{>0}^\times} \\
{\chi_\infty:T(\R) \ra \R^\times} \ar{rr}{\acdot_\infty} &&{\R_{>0}^\times}
\end{tikzcd}
\end{cases}
\]
%%----------------------------------------------------------------------------------------------02
from which an arrow
\[
\begin{cases}
\ \chi_\A:T(\A) \ra \R_{>0}^\times\\
\ \hspace{1.35cm}  x \ra \prod\limits_{p \leq \infty} \chi_p(x_p)
\end{cases}
.
\]

\end{x}
\vspace{0.3cm}

\begin{x}{\small\bf NOTATION} \ %04
\[
T^1(\A) \ = \ \bigcap\limits_{\chi \in X(T)} \Ker \chi_\A.
\]

\end{x}
\vspace{0.3cm}

\begin{x}{\small\bf \un{N.B.}}  \ %05
The  infinite intersection can be replaced by a finite intersection since if 
$\chi_1, \ldots, \chi_d$ is a basis for $X(T)$, then 
\[
T^1(\A) \ = \ \bigcap\limits_{i = 1}^d \ \Ker(\chi_i)_\A.
\]

\end{x}
\vspace{0.3cm}

\begin{x}{\small\bf THEOREM} \ %06
The quotient $T^1(\A) / T(\Q)$ is a compact Hausdorff space.
\end{x}
\vspace{0.3cm}

\begin{x}{\small\bf CONSTRUCTION} \ %07
Let $\Omega_T$ denote the collection of all left invariant $d$-forms on \mT, thus $\Omega_T$ is a 1-dimensional 
vector space over $\Q$.  
Choose a nonzero element $\omega \in \Omega_T$ $-$then $\omega$ determines a left invariant differential form 
of top degree on the $T(\Q_p)$ and $T(\R)$, which in turn determines a Haar measure $\mu_{\Q_p,\omega}$ on the 
$T(\Q_p)$ and a Haar measure $\mu_{\R,\omega}$ on $T(\R)$.

The product
\[
\prod\limits_p \mu_{\Q_p,\omega} (T(\Z_p))
\]
may or may not converge.
\end{x}
\vspace{0.3cm}

%%----------------------------------------------------------------------------------------------03

\begin{x}{\small\bf DEFINITION} \ %08
A sequence $\Lambda = \{\Lambda_p\}$ of positive real numbers is said to be a system of 
\un{convergence coefficients}
\index{convergence coefficients}  
if the product 
\[
\prod\limits_p \ \Lambda_p \mu_{\Q_p,\omega} (T(\Z_p)) 
\]
is convergent.
\end{x}
\vspace{0.3cm}

\begin{x}{\small\bf \un{N.B.}}  \ %09
Covergence coefficients always exist, e.g., 
\[
\Lambda_p \ = \ \frac{1}{\mu_{\Q_p,\omega} (T(\Z_p))}.
\]

\end{x}
\vspace{0.3cm}

\begin{x}{\small\bf LEMMA} \ %10
If the sequence $\Lambda = \{\Lambda_p\}$ is a system of convergence coefficients, then
\[
\mu_{\omega,\Lambda} \ \equiv \ \prod\limits_p \Lambda_p \mu_{\Q_p,\omega} \times \mu_{\R,\omega}
\]
is a Haar measure on $T(\A)$.

\end{x}
\vspace{0.3cm}

\begin{x}{\small\bf \un{N.B.}}  \ %11
Let $\lambda$ be a nonzero rational number $-$then 
\[
\mu_{\Q_p,\lambda\omega} \ = \ \abs{\lambda}_p \mu_{\Q_p,\omega}, 
\quad \mu_{\R,\lambda\omega} \ = \ \abs{\lambda}_\infty \mu_{\R,\omega}.
\]
Therefore
\begin{align*}
\mu_{\lambda\omega,\Lambda} \ 
&\equiv\ 
\prod\limits_p \Lambda_p \mu_{\Q_p,\lambda\omega} \times \mu_{\R,\lambda\omega}
\\[12pt]
&=\ 
\bigl(\hsp \prod\limits_p \abs{\lambda}_p \hsp\bigr)  \hsp
\prod\limits_p \Lambda_p \hsp \mu_{\Q_p,\omega} \times  \abs{\lambda}_\infty \mu_{\R,\omega}
\\[12pt]
&=\ 
\prod\limits_{p \leq \infty} \abs{\lambda}_p \hsp
\prod\limits_p \Lambda_p \hsp \mu_{\Q_p,\omega} \times \mu_{\R,\omega}
\\[12pt]
&=\ \mu_{\omega,\Lambda}.
\end{align*}
%%----------------------------------------------------------------------------------------------04
And this means that the Haar measure $\mu_{\omega,\lambda}$ is independent of the choice of the rational density $\omega$.

Let $\K/\Q$ be a Galois extension relative to which \mT splits $-$then 
\[
X(T_\K) \ = \ \Mor_\K(T_\K,G_{m,\K})
\]
is a $\Gal(\K,\Q)$ lattice.  Call $\Pi$ the representation thereby determined and denote its character by $\chi_{\Pi}$.  
Let
\[
\tL(s,\chi_{\Pi},\K/\Q) \ = \ \prod\limits_p \tL_p(s,\chi_{\Pi},\K/\Q)
\]
be the associated Artin L-function and denote by \mS the set of primes that ramify in $\K$ plus the ``prime at infinity''.
\end{x}
\vspace{0.3cm}

\begin{x}{\small\bf LEMMA} \ %12
$\forall \ p \notin S$, 
\[
\mu_{\Q_p,\omega}(T (\Z_p)) \ = \ \tL_p(1,\chi_{\Pi},\K/\Q)^{-1}.
\]
\end{x}
\vspace{0.3cm}

\begin{x}{\small\bf SCHOLIUM} \ %13
The sequence $\Lambda = \{\Lambda_p\}$ defined by the prescription
\[
\Lambda_p \ = \ \tL_p(1,\chi_\Pi,\K/\Q) \qquad \text{if} \ p \notin S
\]
and 
\[
\Lambda_p \ = \ 1  \qquad \text{if} \ p \in S
\]
is a system of convergence coefficients termed 
\un{canonical}.
\index{convergence coefficients \\ canonical}
\end{x}
\vspace{0.3cm}

\begin{x}{\small\bf LEMMA} \ %14
The Haar measure $\mu_{\omega,\Lambda}$ on $T(\A)$ corresponding to a canonical system of convergence coefficients is 
independent of the choice of $\K$, denote it by $\mu_T$.
\end{x}
\vspace{0.3cm}

\begin{x}{\small\bf DEFINITION} \ %15
$\mu_T$ is the 
\un{Tamagawa measure}
\index{Tamagawa measure}  
on $T(\A)$.
\end{x}
\vspace{0.3cm}

%%----------------------------------------------------------------------------------------------05

Owing to Brauer theory, there is a decomposition of the character $\chi_\Pi$ of $\Pi$ as a finite sum
\[
\chi_{\Pi} \ = \ d_{\chi_0} + \sum\limits_{j = 1}^M m_j \chi_j,
\]
where $\chi_0$ is the principal character of $\Gal(\K/\Q)$ $(\chi_0(\sigma) = 1$ for all $\sigma \in \Gal(\K/\Q))$, 
the $m_j$ are positive integers, and the $\chi_j$ are irreducible characters of $Gal(\K/\Q)$.  
Standard properties of Artin L-functions then imply that 
\[
\tL(s,\chi_{\Pi},\K/\Q) \ = \ \zeta(s)^d \ \prod\limits_{j = 1}^M \ \tL(s, \chi_j,\K/\Q)^{m_j}.
\]
\vspace{0.3cm}

\begin{x}{\small\bf FACT} \ %16
\[
\ \tL(1, \chi_j,\K/\Q)^{m_j} \ \neq \ 0 \qquad (1 \leq j \leq M).
\]
Therefore
\begin{align*}
\lim\limits_{s \ra 1} (s - 1)^d \tL(s,\chi_{\Pi},\K/\Q) \ 
&= \ \prod\limits_{j = 1}^M \ \tL(1, \chi_j,\K/\Q)^{m_j}\\
&\neq \ 0.
\end{align*}
\end{x}
\vspace{0.3cm}

\begin{x}{\small\bf LEMMA} \ %17
The limit on the left is positive and independent of the choice of $\K$, denote it by $\rho_T$.
\end{x}
\vspace{0.3cm}

\begin{x}{\small\bf DEFINITION} \ %18
$\rho_T$ is the 
\un{residue}
\index{residue} 
of \mT.

\vspace{0.3cm}

Define a map
\[
T:T(\A) \ \ra \ (\R_{>0}^\times)^d
\]
%%----------------------------------------------------------------------------------------------06
by the rule
\[
T(x) \ = \ ((\chi_1)_\A(x), \ldots, (\chi_d)_\A(x)).
\]
Then the kernel of \mT is $T^1(\A)$, hence \mT drops to an isomorphism
\[
T^1:T(\A)/T^1(\A) \ \ra \ (\R_{>0}^\times)^d.
\]

\end{x}
\vspace{0.3cm}

\begin{x}{\small\bf DEFINITION} \ %19
The 
\un{standard measure}
\index{standard measure} 
on $T(\A)/T^1(\A)$ is the pullback via $T^1$ of the product measure
\[
\prod\limits_{i = 1}^d \ \frac{d_{t_i}}{t_i}
\]
on $(\R_{>0}^\times)^d$.
\vspace{0.3cm}

Consider now the formalism
\[
d(T(\A)) \ = \ d(T(\A) / T^1(\A)) d(T^1(\A) / T(\Q)) d(T(\Q))
\]
in which:
\vspace{0.3cm}

\qquad\qquad \textbullet \quad $d(T(\A))$ is the Tamagawa measure on $T(\A)$ multiplied by $\ds\frac{1}{\rho_T}$.
\vspace{0.2cm}

\qquad\qquad \textbullet \quad $d(T(\A) / T^1(\A))$ is the standard measure on $T(\A) / T^1(\A)$.
\vspace{0.2cm}

\qquad\qquad \textbullet \quad $d(T(\Q))$ is the counting measure on $T(\Q)$.
\end{x}
\vspace{0.3cm}

\begin{x}{\small\bf DEFINITION} \ %20
The 
\un{Tamagawa number}
\index{Tamagawa number} 
$\tau(T)$ is the volume
\[
\tau(T) \ = \ \int\limits_{T^1(\A)/T(\Q)} 1
\]
of the compact Hausdorff space $T^1(\A)/T(\Q)$ per the invariant measure
\[
d(T^1(\A)/T(\Q))
\]
%%----------------------------------------------------------------------------------------------07
such that
\[
\frac{\mu_T}{\rho_T} \ = \ d(T(\A) / T^1(\A)) \hsp d(T^1(\A)/T(\Q)) \hsp d(T(\Q)). %dmc error orig ? missing /
\]

\end{x}
\vspace{0.3cm}

\begin{x}{\small\bf \un{N.B.}}  \ %21
To be completely precise, the integral formula
\[
\int\limits_{T(\A)} \ = \ \int\limits_{T(\A)/T^1(\A)} \int\limits_{T^1(\A)}
\]
fixes the invariant measure on $T^1(\A)$ and from there the integral formula
\[
 \int\limits_{T^1(\A)} \ = \  \int\limits_{T^1(\A)/T(\Q)} \int\limits_{T(\Q)}
\]
fixes the invariant measure on $T^1(\A)/T(\Q)$ , its volume then being the Tamagawa number $\tau(T)$.

\vspace{0.3cm}

[Note: \ 
If \mT is $\Q$-anisotropic, then $T(\A) = T^1(\A)$.]
\end{x}
\vspace{0.2cm}

\begin{x}{\small\bf EXAMPLE} \ %22
Take $T = G_{m,\Q}$ and $\omega = \ds\frac{dx}{x}$ $-$then 
\[
\vol_{\frac{dx}{\abs{x}_p}} (\Z_p^\times) \ = \ \frac{p-1}{p} \ = \ 1 - \frac{1}{p}
\]
and the canonical covergence coefficients are the
\[
\bigl(1 - \frac{1}{p}\bigr)^{-1}.
\]
Here $d = 1$ and 
\[
\lim\limits_{s \ra 1} \ (s - 1) \zeta(s) \ = \ 1 \ \implies \ \rho_T = 1.
\]
Working through the definitions, one concludes that $\tau(T) = 1$ or still, 
\[
\vol(\I^1/\Q^\times) \ = \ 1.
\]
\end{x}
\vspace{0.3cm}

%%----------------------------------------------------------------------------------------------08

\begin{x}{\small\bf REMARK} \ %23
Take $T = \Res_{\K/\Q}(G_{m,\K})$ $-$then it turns out that $\tau(T)$ is the Tamagawa number of $G_{m,\K}$ 
computed relative to $\K$ (and not relative to $\Q$ \ldots).  
From this, it follows that $\tau(T) = 1$, matters hinging on the ``famous formula''
\[
\lim\limits_{s \ra 1} \ (s - 1) \zeta_\K(s) \ = \ 
\frac{2^{r_1}(2\pi)^{r_2}}{w_\K \abs{d_\K}^{1/2}} h_\K R_\K.
\]
\end{x}
\vspace{0.3cm}

\begin{x}{\small\bf LEMMA} \ %24
Let \mF be an integrable function on $(\R_{>0}^\times)^d$ $-$then 
\[
\tau(T) \ = \  
\frac{\ds{\frac{1}{\rho_T}} \hsp \int\limits_{T(\A)/T(\Q)} F(T(x)) d\mu_T(x)}
{\ds\int\limits_{(\R_{>0}^\times)^d} F(t_1, \ldots, t_d) \frac{dt_1}{t_1}\ldots\frac{dt_d}{t_d}}\hsp.
\]
\end{x}
\vspace{0.3cm}

\begin{x}{\small\bf EXAMPLE} \ %25
Let $T = G_{m,\Q}$ $-$then
\[
\tau(T) \ = \ \ds\frac{\ds\int\limits_{\I/\Q^\times} F(\abs{x}_\A) d\mu_T(x)}{\ds\int_0^\infty \frac{F(t)}{t}dt},
\]
$\rho_T$ being 1 in this case.  
To see that $\tau(T) = 1$, make the calculation by choosing 
\[
F(t) \ = \ 2 t e^{-\pi t^2}.
\]
\vspace{0.2cm}

[Note: \ 
Recall that 
\[
\prod\limits_p \Z_p^\times \times \R_{>0}^\times
\]
is a fundamental domain for $\I/\Q^\times$.]
\end{x}
\vspace{0.3cm}

\begin{x}{\small\bf NOTATION} \ %26
Put
\[
H^1(\Q,T) \ = \ H^1(\Gal(\Q^\sep/\Q), T(\Q^\sep))
\]
%%----------------------------------------------------------------------------------------------09
and for $p \leq \infty$,
\[
H^1(\Q_p,T) \ = \ H^1(\Gal(\Q_p^\sep/\Q_p), T(\Q_p^\sep)).
\]
\end{x}
\vspace{0.3cm}

\begin{x}{\small\bf LEMMA} \ %27
There is a canonical arrow 
\[
H^1(\Q,T) \ \ra \ H^1(\Q_p,T).
\]

PROOF \ 
Put
\[
G \ = \ \Gal(\ov{\Q}/\Q) \qquad (\ov{\Q} = \Q^\sep)
\]
and
\[
G_p \ = \ \Gal(\ov{\Q}_p/\Q_p) \qquad (\ov{\Q}_p = \Q_p^\sep).
\]
Then schematically
\[
\begin{tikzcd}[sep=large]
&{\ov{\Q}_p} 
\arrow[ld,dash] \arrow[rd,dash,"G_p"]\\
{\ov{\Q}} &&{\Q_p} \\
&{\Q} \arrow[lu,dash,"G"] \arrow[ru,dash] 
\end{tikzcd}
.
\]
\vspace{0.2cm}

\qquad\qquad 1. \quad 
There is an arrow of restriction
\[
\rho:G_p \ \ra \ G
\]
and a morphism $T(\Q) \ra T(\ov{\Q}_p)$ of $G_p$-modules, $T(\Q)$ being viewed as a $G_p$-module via $\rho$.
\vspace{0.2cm}

\qquad\qquad 2. \quad 
The canonical arrow
\[
H^1(\Q,T) \ \ra \ H^1(\Q_p,T)
\]
is then the result of composing the map
\[
H^1(G,T(\Q)) \ \ra \ H^1(G_p,T(\Q))
\]

%%----------------------------------------------------------------------------------------------10
with the map
\[
H^1(G_p,T(\Q)) \ \ra \ H^1(G_p,T(\ov{Q}_p)).
\]
\end{x}
\vspace{0.3cm}

\begin{x}{\small\bf NOTATION} \ %28
Put
\[
\Sh(T) \ = \ \Ker(H^1(\Q,T) \ \ra \ \coprod\limits_{p \leq \infty} H^1(\Q_p,T)).
\]
\end{x}
\vspace{0.3cm}

\begin{x}{\small\bf DEFINITION} \ %29
$\Sh(T)$ is the 
\un{Tate-Shafarevich group}
\index{Tate-Shafarevich group} 
of \mT.
\end{x}
\vspace{0.3cm}

\begin{x}{\small\bf THEOREM} \ %30
$\Sh(T)$ is a finite group.
\end{x}
\vspace{0.3cm}

\begin{x}{\small\bf EXAMPLE} \ %31
If $\K$ is a finite extension of $\Q$, then 
\[
H^1(\Q,\Res_{\K/\Q}(G_{m,\K})) \ = \ 1.
\]
Therefore in this case
\[
\#(\Sh(T)) \ = \ 1.
\]
\end{x}
\vspace{0.3cm}

\begin{x}{\small\bf REMARK} \ %32
By comparison,
\[
H^1(\Q,\Res_{\K/\Q}^{(1)}(G_{m,\K})) \ \approx \ \Q^\times / N_{\K/\Q} (\K^\times).
\]

[Consider the short exact sequence
\[
\begin{tikzcd}%[sep=large]
{1} \ar{r} 
&{\Res_{\K/\Q}^{(1)}(G_{m,\K})} \ar{r} 
&{\Res_{\K/\Q}(G_{m,\K})} \ar{rr}{\tN_{\K/\Q}} 
&&{G_{m,\Q}} \ar{r} 
&{1}
\end{tikzcd}
.]
\]

\end{x}
\vspace{0.3cm}

\begin{x}{\small\bf NOTATION} \ %33
Put
\[
\LY(T) \ = \ \CoKer \bigl(H^1(\Q,T) \ \ra \ \coprod\limits_{p \leq \infty} H^1(\Q_p,T)\bigr).
\]
\end{x}
\vspace{0.3cm}

%%----------------------------------------------------------------------------------------------11

\begin{x}{\small\bf THEOREM} \ %34
$\LY(T)$ is a finite group.
\end{x}
\vspace{0.3cm}

\begin{x}{\small\bf MAIN THEOREM} \ %35
The Tamagawa number $\tau(T)$ is given by the formula
\[
\tau(T) \ = \ \frac{\# (\LY(T))}{\# (\Sh(T))}.
\]
\end{x}
\vspace{0.3cm}

\begin{x}{\small\bf EXAMPLE} \ %36
If $\K$ is a finite extension of $\Q$, then 
\[
H^1(\Q_p,\Res_{\K/\Q}(G_{m,\K})) \ = \ 1.
\]
Therefore in this case
\[
\# (\LY(T)) \ = \ 1.
\]
\vspace{0.3cm}

It follows from the main theorem that $\tau(T)$ is a positive rational number.  
Still, there are examples of finite abelian extensions $\K / \Q$ such that 
\[
\tau(\Res_{\K/\Q}^{(1)} G_{m,\K})
\]
is not a positive integer.
\end{x}
\vspace{0.3cm}

%%----------------------------------------------------------------------------------------------12

%%%%%%%%%%%%%%%%%%%%%%%%%%%%%%%%%%%%%%
%%%%%%%%%%%%%%%%%%%%%%%%%%%%%%%%%%%%%%
%%%%%%%%%%%%%%%%%%%%%%%%%%%%%%%%%%%%%%

%% file: _references.tex
\renewcommand{\thepage}{References-\arabic{page}}
\setcounter{page}{1}
\begingroup
\center {\textbf{REFERENCES}}\\
\endgroup
\vspace{0.5cm}
%\[
%\textbf{References}
%\]

\noindent G\"ortz, U. et al.

[1] \quad 
Algebraic Geometry I., VieWeg + Teubner Verlag, 2010.
\vspace{0.5cm}
%%%%%%%%%%%%%%%%%%%

\noindent Langlands, R. P. 

[2] \quad
The Representation Theory of Abelian Algebraic Groups, Pacific Journal of Mathematics, vol. 181, Issue 3, 1997, pp. 231-250.
\vspace{0.5cm}
%%%%%%%%%%%%%%%%%%%

\noindent Milne, J. S.

[3] \quad
Algebraic Groups: The Theory of Group Schemes of Finite Type over a Field, Cambridge University Press, 2017.
\vspace{0.5cm}
%%%%%%%%%%%%%%%%%%%

\noindent Ono, T.

[4-(a)] \quad
Arithmetic of Algebraic Tori, Annals of Mathematics, vol. 74, 1961, pp. 101-139.
\vspace{0.1cm}

[4-(b)] \quad
On the Tamagawa Number of Algebraic Tori, Annals of Mathematics, vol. 78, 1963, pp. 47-73.
%%%%%%%%%%%%%%%%%%%